\documentclass[12pt,twoside]{article} 
\usepackage[utf8]{inputenc}
\usepackage{graphicx}
\usepackage{enumitem}
\usepackage{latexsym}
\usepackage{amssymb}
\usepackage{epsfig}
\usepackage{xcolor}
\usepackage{amsmath}
\usepackage{float}
\usepackage{caption}
\usepackage{xr}

\makeatletter

\@addtoreset{equation}{section}
\makeatother

\allowdisplaybreaks 

\def\en{\mathbb{N}}

\def\er{\mathbb{R}}

\def\P{\mathbb{P}}
\def\E{\mathbb{E}}

\def\e{\varepsilon}

\def\boxi{\begin{flushright}$\Box$\end{flushright}}
\def\beq{\begin{eqnarray*}}
\def\eeq{\end{eqnarray*}}

\newcommand{\nto}{\xrightarrow[n\to\infty]{}}
\def\bx{{\bf x}}
\def\bX{{\bf X}}
\def\bt{{\bf t}}
\def\bj{{\bf j}}
\def\bl{{\bf l}}

\newtheorem{theo}{Theorem}[section]
\newtheorem{lemma}[theo]{Lemma}
\newtheorem{cor}[theo]{Corollary}
\newtheorem{rem}[theo]{Remark}

\newtheorem{proposition}[theo]{Proposition}

\newtheorem{remark}[theo]{Remark}
\newtheorem{theorem}[theo]{Theorem}

\begin{document}

\title{\bf Multivariate boundary regression models}

\author{{\sc Leonie Selk$^1$}, {\sc Charles Tillier$^2$} and {\sc Orlando Marigliano$^3$}\\{\small $^1$Department of Mathematics, University of Hamburg} \\{\small $^2$Laboratoire de Mathématiques de Versailles, University of Versailles-Saint-Quentin-en-Yvelines}\\ {\small $^3$Max Planck Institute for Mathematics in the Sciences}}


\maketitle

	\begin{abstract}
	
	In this work, we consider  a  multivariate regression model with  one-sided  errors. We assume for the regression function to lie in a general H\"{o}lder class and estimate it via a nonparametric local polynomial approach that consists of minimization of the local integral of a polynomial approximation lying above the data points.
	While the consideration of multivariate covariates offers an undeniable opportunity from an application-oriented standpoint, it requires a new method of proof to replace the established ones for the univariate case. 
	
	The main purpose of this paper is to show the uniform consistency and to provide the rates of convergence of the considered nonparametric estimator for both multivariate random covariates and multivariate deterministic design points. To demonstrate the performance of the estimators, the small sample behavior is investigated in a simulation study in dimension two and three.

	
\end{abstract}


\section{Introduction}
We consider nonparametric regression models with one-sided errors that take the general form
\begin{align}\label{intro:def_model}
Y_i=g(X_i)+\varepsilon_i, \ \ \ i=1, \ldots,n
\end{align}
where $Y_i$ is the response variable, $X_i$ is the multivariate random or deterministic covariate, $g$ is the  unknown regression function corresponding to the upper boundary curve and $\varepsilon_i$ is a nonpositive random error term. The statistical issue of such \textit{boundary regression models} (BRM) lies on the frontier estimation, in other words, on the estimation of $g$ based on the observations $(X_i,Y_i)$, $i=1, \ldots,n$ where $n$ is the sample size of the available data.

BRM  have received an increasing attention in the last past years and are closely related to production frontier models (PFM): both share the objective of the estimation of the frontier - the boundary curve - and contribute to the same applications. PFM appeared in the seminal contribution of Farrell (1957) in  data envelopment analysis (DEA) to answer the need of developing nonparametric methods to assess efficiency (i.e.\ production of the maximum output which is  feasible for the given inputs) of a system. DEA and its extensions developed e.g.\ in Deprins, Simar and Tulkens (1984) and Lovell et al.\ (1994) are now recognized as powerful tools for evaluating the performance of a system and have countless applications, among others, in social sciences, health-care evaluation systems and banking sectors. We refer to the books Cooper, Seiford and Zhu  (2011) and Ramanathan (2003) for a comprehensive treatment of such methods and an exhaustive development of the applications. On the other hand, stochastic frontier analysis formulated independently by Aigner, Lovell and Schmidt (1977) and Meeusen and van Den Broeck (1977) offers an interesting alternative with parametric estimations of the frontier; see also the books of Kumbhakar and  Lovell (2003) and Cornwell and Schmidt (2008) for more recent references.

There is vast literature on PFM dealing with the estimation of the boundary curve. Numerous parametric and non parametric techniques have been proposed, for instance using extreme-value based estimators (Girard and Jacob (2003), de Haan and Resnick (1994), Hall, Nussbaum and Stern (1997), Menneteau (2008), Girard and Jacob (2004),  Gardes (2002) and Gijbels and Peng (2000)), projections techniques (Jacob and Suquet (1995)), kernel based estimators (Girard, Guillou and Stupfler (2013), Girard and Jacob (2008)), maximum likelihood estimators (Kumbhakar et al.\ (2007) and Simar and Zelenyuk (2011)). Some estimators need the boundary curve to be monotone; see e.g.\ Daouia and Simar (2005), Daouia, Noh and Park (2016) and Gijbels et al.\ (1999).

In contrast to the aforementioned methods, we concentrate in this paper on an alternative approach that consists in  approximating the regression function $g$ locally by a polynomial lying above the data points.  Polynomial estimators in frontier estimation have been widely studied in the literature and are employed in several works, for instance in Hall, Park and Stern (1998), Hall and Park (2004), Girard, Iouditski and Nazin (2005), Knight (2001) and in Hall and Van Keilegom (2009); see also the literature using the alternative method of piecewise polynomials (e.g.\  Korostelëv, Simar and Tsybakov (1995), H\"ardle, Park and Tsybakov (1995)  and Chapter 3 in Korostelëv and Tsybakov (1993), as well as Tsybakov (1994) and Chapter 5  in Korostelëv and Tsybakov (1993) for the multivariate case) in which the estimation of the boundary of sets of the form $\{(x,y):0\leq x\leq 1, 0\leq y\leq g(x)\}$ is considered.

In the context of BRM, local polynomial estimator of the frontier benefit from attractive properties. Under the assumptions of  H\"older boundary curve and regularly varying errors, Jirak, Meister and Reiß (2014) suggested an adaptive estimator for the  boundary curve $g$ using a local polynomial estimation based on local extreme value statistics. An adaptive procedure - a fully data-driven estimation procedure - is constructed by applying a nested version of Lepski method which shows no loss in the convergence rates with respect to the general $L_q$-risk. Drees, Neumeyer and Selk (2019) estimated the regression function similarly to Jirak et al.\ (2014) via minimization of the local integral of a polynomial approximation in the context of equidistant design points. By showing uniform rates of convergence for the  regression estimator they proposed distribution-free tests  of error distributions where the test statistics are based on empirical processes of residuals. They also discussed asymptotically distribution-free hypotheses tests for independence of the error distribution from the points of measurement and for monotonicity of the boundary function as well. Relatedly, Neumeyer, Selk and Tillier (2019) showed consistency of such estimators under mild condition.

However, while frontier estimators based on local polynomial approximations appears to be an interesting  method, this procedure suffers from the fact the literature studies have only been carried out for univariate samples. For real data application, such a setup of a single covariate means that the output has to be explained by a single input which is not realistic and hence the current theoretical guarantees do not apply. For instance, in medicine, Narimatsu et al.\ (2015) focus on predicting obesity investigating the inverse of body mass index (the output) using two inputs, the number of calories expended and ingested. In agriculture, Kelly et al.\ (2012) studied the efficiency of dairy farms thanks to several inputs: the land size, the cow numbers, the labour and different costs.  More generally, in PFM, no matter what application field one is interested in, all the studies are conducted with several inputs. We refer the reader to the preceding paragraph dedicated to DEA for further applications and to the references within.  As a consequence, until now this technique does not bring any concrete method from an application point of view and from this perspective, extending to the multivariate setup is important.

On the other hand, beyond the practical interest this generalization raises, investigating the asymptotic properties to higher dimensions is also interesting from a theoretical standpoint. Indeed, Jirak et al.\ (2014)  and Dress et al.\ (2019) obtain asymptotic results on the univariate case by passing to the complex numbers and using the fundamental theorem of algebra. Since this technique is not available to us in the multivariate case, the extension is quite demanding and requires a  new approach. Besides, the proposed approach in this paper provides new proofs for the univariate case and the results of Jirak et al.\ (2014)  and Dress et al.\ (2019) follow as special cases.  Further details are given in Remark \ref{rem:extension}.

More often than not, deterministic covariates and especially  equidistant fixed design points are considered. In contrast in the paper at hand, we investigate the two cases of random and deterministic covariates which are both of particular interest. Deterministic covariates are often used in real-life applications when time is involved in the data set.  For instance Jirak et al.\ (2014) studied the monthly sunspot observations and the annual best running times of 1500 meters; see also the plentiful applications in energy and environmental performance analysis provided in Mardani et al.\ (2018). Besides, deterministic design is met across a number of papers in regression models, see for instance Brown and Low (1996), Meister and Rei\ss \ (2013) and the references within. The case of random covariates is obviously the most relevant and appears in essence in many applications, among other, in insurance and finance when analyzing optimality of portfolios and efficiency strategies; see also the extensive literature on modern portofolio theory (e.g.\ Francis and Kim (2013) and Goetzmann et al.\ (2014)).

In light of these motivations, the main aim of this paper is to show uniform consistency and to provide rates of convergence of an estimator based on the minimization of the local integral of a polynomial lying above the data points  for both \textit{multivariate random covariates} and \textit{multivariate deterministic design points}, under the main assumptions of regular variation of the nonpositive errors and $\beta$-H\"older class of the boundary curve. 


The remaining part of the manuscript is organized as follows. In section \ref{model} the model is explained, while in section \ref{estimators} the estimation procedure is described. In section \ref{main} we show uniform consistency and provide rates of convergence of the estimator of the regression function for both random and deterministic multivariate covariates. Section \ref{simulations} is dedicated to a simulation study to investigate the small sample behavior in dimension two and three. The proofs are summed up in section \ref{appendix}. 

\subsection*{Notation and shortcuts}

\underline{Notation}:  $\mathcal{B}$ stands for Borelian sets; $\lfloor \cdot \rfloor$ and $\lceil \cdot \rceil$ are the floor and ceiling functions respectively; $\langle x\rangle$ means the largest natural number that is strictly smaller than $x$; $\bar{F}=1-F$ denotes the survival function associated to a cdf $F$;  $\Vert \cdot \Vert_{\infty}$ denotes the supremum norm; $X_1 \overset{d}{=} X_2$ means that two random variables $X_1,X_2$ share the same distribution; $a_n \underset{n \to \infty}{\sim} b_n$ holds if $\lim_{n \to \infty}a_n/b_n=1$ for two sequences $(a_n)_{n \geq 1}$ and $(b_n)_{n \geq 1}$ of nonnegative numbers. Generally, vectors are highlighted in bold writing. For vectors $\bx\in\er^q$ let $x^{(r)}$ denote the $r$-th component of $\bx$ for $r=1,\ldots,q$. By $\|\bx\|$ we mean the maximum norm that is $\|\bx\|:=\max_{r\in\{1,\ldots,q\}}|x^{(r)}|$. 
For multivariate polynomials we use the multiindex notation where for a vector $\bx\in\er^q$ and a multiindex $\bj=(j_1,\ldots,j_q)\in\en_0^q$ we define $|\bj|:=j_1+\ldots+j_q$, $\bj!:=j_1!\cdot\ldots\cdot j_q!$ and $\bx^\bj:={x^{(1)}}^{j_1}\cdot\ldots\cdot{x^{(q)}}^{j_q}$.
	
\underline{Shortcuts}: \textit{cdf} stands for cumulative distribution function and \textit{iid} for independent and identically distributed.
	
	\section{The model}\label{model}
	We focus on  nonparametric boundary regression models with univariate observations $Y_i$ and  multivariate covariates $X_i$ of the form
\begin{align}\label{Model:def_model}
Y_i=g(X_i)+\varepsilon_i, \ \ \ i=1, \ldots,n
\end{align}
where the errors $\varepsilon_i$ are iid non-positive univariate random variables that are independent of the $X_i$ and $n$ stands for the sample size. The unknown regression function $g$ thus corresponds to the upper boundary curve.

Independence of errors $\varepsilon_i$ and covariates $X_i$ is a typical assumption in regression models and is met among others in Müller and Wefelmeyer (2010), Meister and Rei{\ss} (2013), Rei{\ss} and Selk (2017) and Drees et al.\ (2019). Such assumption is crucial  and often needed in the framework of frontier models when using statistical methods such as bootstrap procedures; see  Wilson (2003) and  Simar and Wilson (1998). In the case this hypothesis does not hold, one may consider parametric transformations e.g.\ exponential, Box-Cox (Box and Cox (1964)), or  sinh-arcsinh transformations (Jones and Pewsey (2009)) of the response variable $Y_i$ in order to retreive the framework of independence between errors and  covariates; see  Neumeyer et al.\ (2019) for the univariate case and also Linton,  Sperlich and Van Keilegom  (2008).

In the paper at hand, we investigate both cases of random and fixed design points. The former  means that the covariates $X_i$ are $q$-dimensional random variables while the latter assumes that $X_i$ are  deterministically spread over $[0,1]^q$. Without loss of generality, we work for ease of reading with design points lying on $[0,1]^q$ but results extend effortless to any cartesian product  of  one-dimensional closed intervals. 
		
	\subsection{The random design case}

In the random design case we consider the nonparametric boundary regression model with  independent and identically distributed  observations $(\bX_i,Y_i)$, $i=1,\ldots,n$  defined by
\begin{equation}\label{model-random}
Y_i=g(\bX_i)+\e_i
\end{equation}
corresponding to model \eqref{Model:def_model}, where the design points $\bX_i$ are multivariate random covariates distributed on $[0,1]^q$ that fulfill assumption (K4). The errors $\e_i$ are assumed to be iid non-positive random variables that satisfy (K2). The precise statements of assumptions (K2) and (K4) are given in Section 3.1.

\subsection{The fixed design case}\label{model-fixedequi}

In the fixed design case we consider a triangular array of independent observations $Y_{i,n}$
and deterministic design points $\bx_{i,n}$ in $[0,1]^q$ for $i=1,\ldots,n$. Thus we conduct the nonparametric boundary regression model
\begin{equation}\label{model0}
Y_{i,n}=g(\bx_{i,n})+\varepsilon_{i,n},
\end{equation}
corresponding to model \eqref{Model:def_model} with errors $\varepsilon_{i,n}$ that are univariate non-positive independent and identically distributed random variables that satisfy assumption (K2).

We allow for fixed equidistant as well as fixed nonequidistant design. In the first case we consider $\bx_{1,n},\ldots,\bx_{n,n}$ that form a grid 
\[\begin{pmatrix}n^{-\frac 1q}\\ n^{-\frac 1q}\\ \vdots\\n^{-\frac 1q}\end{pmatrix},\begin{pmatrix}2n^{-\frac 1q}\\ n^{-\frac 1q}\\\vdots\\n^{-\frac 1q}\end{pmatrix},\ldots,\begin{pmatrix}1\\\vdots\\1\\(n^{\frac 1q}-1)n^{-\frac 1q}\end{pmatrix}, \begin{pmatrix}1\\\vdots\\1\\1\end{pmatrix}\]
where we assume that $n^{\frac 1q}$ is an integer.  Note that when $q=1$ the univariate equidistant design simplifies to $\bx_{i,n}= i/n$ for $i=1, \ldots,n$. In the second case the points are not necessarily equidistant, but we assume that they are even enough distributed on $[0,1]^q$, see Assumption (K4') below.

\section{Estimating the regression function}\label{estimators}
To estimate the boundary curve, we use an estimator that locally approximates the regression function $g$ by a polynomial lying above the data points; see Theorem 2.2 in Drees et al.\ (2019) for further details in the univariate case.

\subsection{The random design case}
For $\bx \in [0,1]^q$, we consider the regression function estimator $\hat{g}$ defined as  
\begin{align}\label{estimatorIntegral}
\hat{g}_n(\bx):=\hat{g}(\bx):=p(\bx)
\end{align}
where $p$ is a multivariate polynomial of total degree $\beta^* \in \mathbb{N}_0$ and minimizes the local integral 
\begin{align}\label{eq:integ}
\int_{[0,1]^q} p(\bt)I\{\|\bt-\bx\|\leq h_n\}d\bt
\end{align}
under the constraints $p(\bX_i) \geq Y_{i}$ for $\|\bX_i-\bx\|\leq h_n$.
Here, $h_n > 0$ is taken to satisfy assumption (K3) below.\newline

\begin{rem}\label{Rem-solvePoly}
The polynomial $p$ is the solution to the linear optimization problem
\begin{align}\label{optProb}
    \text{\rm minimize} &\quad v^T  p \\*
    \text{\rm subject to} &\quad A p \geq y, \nonumber
\end{align}
where $p$ is represented by its vector of coefficients, $v$ is the vector representing the linear functional $\int_{||\bt-\bx||\leq h_n}$, the matrix $A$ is the multivariate Vandermonde matrix whose $i$-th row has as its entries all the monomials of degree at most $\beta^*$ in the entries of $\bX_i$, and $y$ is the vector with $y_i = Y_i$. For the estimator $\hat g$ to be well-defined, it is necessary that this problem is bounded from below, that is, that the objective is bounded from below on the polytope defined by the constraints. This need not always be the case, as the example in Jirak et al.\ (2014) with $\beta^*=2$ and two support points demonstrates. However, the alternate optimization problem proposed in Jirak et al.\ (2014) has the same problem of unboundedness. When $q=1$, Problem \eqref{optProb} is bounded whenever we have at least $\beta^*+1$ points. This follows from the fact that the univariate Vandermonde matrix is totally positive when the support points $X_i$ are positive. However, for higher dimensions it is not as simple. There, the Vandermonde matrix with
$\binom{q+\beta^*}{q}$
rows needs not be invertible. As of now, we believe that for $q>1$ the boundedness of \eqref{optProb} needs to be checked on a case-by-case basis using linear optimization algorithms. By duality theory we know that \eqref{optProb} is bounded if and only if there exists a vector $g\geq 0$ with $A^Tg = v$. This linear program can become very large as $q$ and $\beta^*$ grow. For instance, if $q=\beta^*$ and $A$ has $N = \binom{2q}{q}$ rows, then the interior point algorithm from Vaidya (1989) runs in $O(N^{2.5})$ time in the worst case, which in terms of $q$ grows faster than $4^{2.5 q}\cdot q^{-1.25}$. Nevertheless, this is a theoretical worst case and the average case might be better, possibly using another algorithm. For implementations, the authors suggest using the Python module \texttt{scipy.optimize.linprog}, which by default uses an interior point algorithm based on the MOSEK interior point optimizer by Andersen and Andersen (2000).
\end{rem}

\begin{rem}
To illustrate the estimation procedure we take a look at the simplest case of $\beta^*=0$ which results in a local constant approximation. The estimator defined in \eqref{estimatorIntegral}-\eqref{eq:integ} then simplifies to $\hat g(\bx)=\max\{Y_{i}\, |\, i=1,\ldots,n\ \text{with}\ \|\bX_i-\bx\|\leq h_n\}$. In Figure \ref{figlocConst} an example for $q=1$ and uniformly distributed $\bX_i$ is shown. For each $\bx$-value a constant function is fitted to the data in the neighborhood.
\begin{figure}
\begin{center}
\begin{minipage}[l]{0.35\textwidth}
\includegraphics[width=\textwidth]{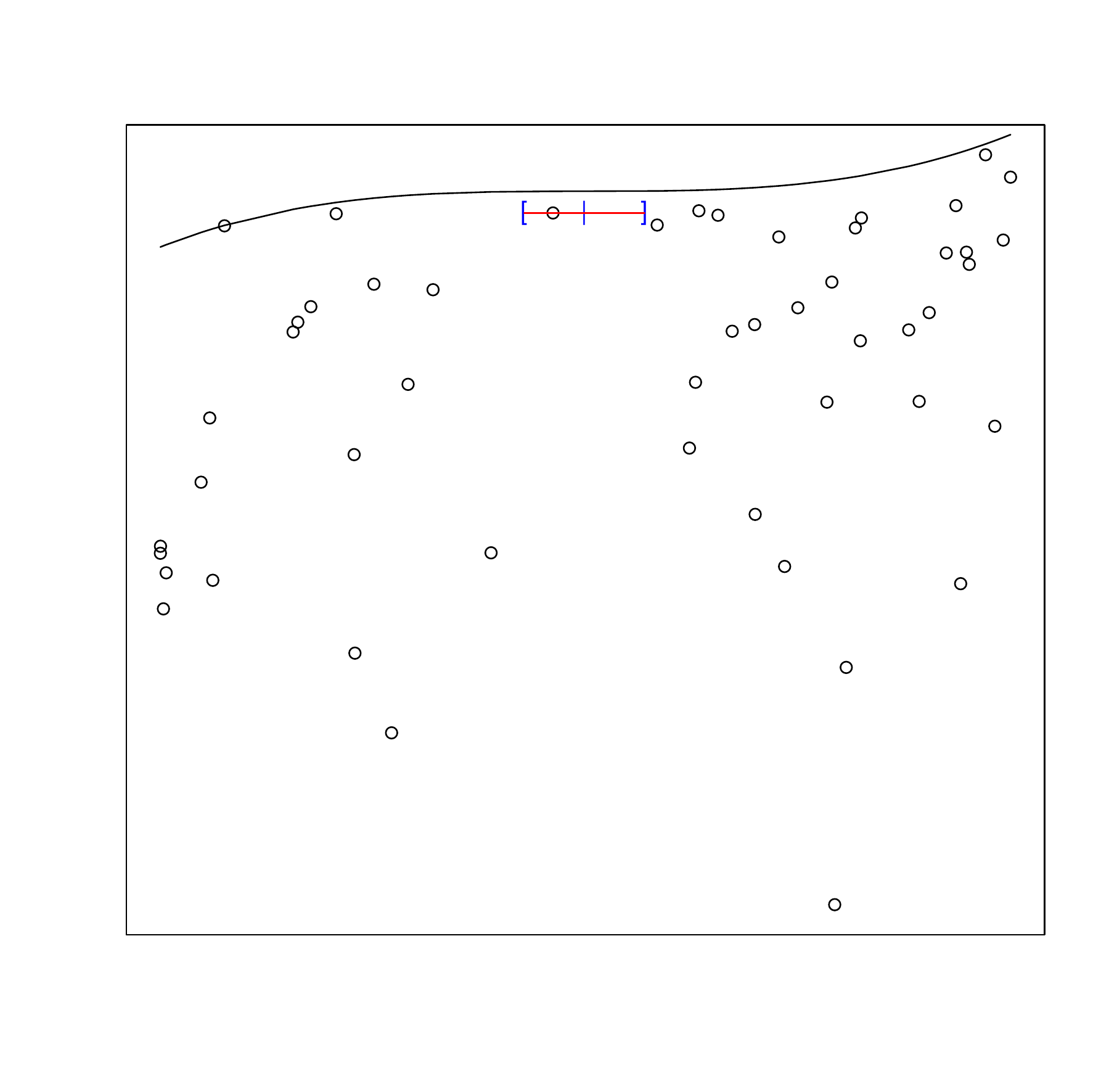}
\end{minipage}
\begin{minipage}[r]{0.35\textwidth}
\includegraphics[width=\textwidth]{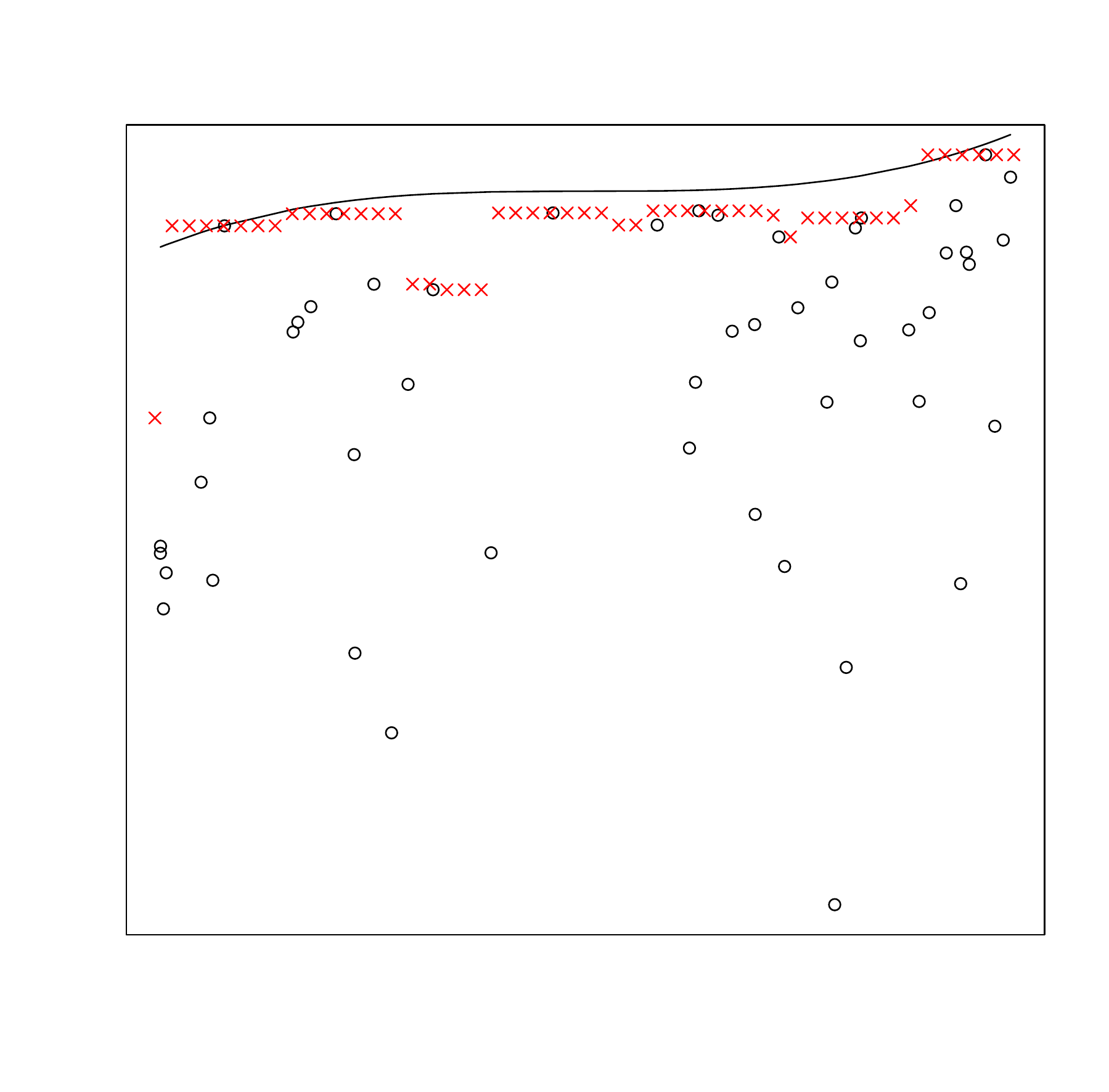}
\end{minipage}
\caption{Scatter plot of $(\bX_i,Y_i)$, $i=1,\ldots,n$ together with the true regression function (black solid curve) in $q=1$. A local constant approximation is considered. On the left hand side the blue vertical bar marks the point $(\bx,\hat g(\bx))$ for some given $\bx$ and the blue square brackets the endpoints of the interval $[\bx-h_n,\bx+h_n]$. On the right hand side the red crosses indicate the estimated values $\hat g(\bx)$ for different values of $\bx$.}\label{figlocConst}
\end{center}
\end{figure}
\end{rem}

We work under the following four assumptions (K1)-(K4).
\begin{enumerate}
\item[(K1)] \textbf{Regression function}: $g$ belongs to some H\"older class of order $\beta \in (0, \infty)$ that is $g$ is $\langle \beta \rangle$-times differentiable on $[0,1]^q$ and all partial derivatives of order $\langle \beta \rangle$ satisfy
	\begin{align}\label{def:holder}
	|D^\bj g(\bt)-D^\bj g(\bx)|\leq c_g\|\bt-\bx\|^{\beta-\langle \beta \rangle}\qquad\forall \bx,\bt\in [0,1]^q, \ \forall\bj\in\en_0^q \text{ with } |\bj|=\langle \beta\rangle
	\end{align}
for some $c_g<\infty$ where
\[D^\bj=\frac{\partial^{j_1+\ldots+j_q}}{\partial x_1^{j_1}\ldots\partial x_q^{j_q}}.\]
		
\item[(K2)] \textbf{Errors distribution}: The errors $\varepsilon_i$ are independent and identically distributed  on $(-\infty,0]$ with common cdf  $F$ that satisfies
	\begin{align*}
	\overline{F}(y)=c|y|^\alpha + r(y), \ \ \ y<0,
	\end{align*} 
	with $\alpha,c>0$ and $r(y)=o(|y|^\alpha)$ when $y\nearrow 0$.
		
\item[(K3)] \textbf{Bandwidths}:  $(h_n)_{n\in\mathbb{N}}$ is a sequence of positive bandwidths that satisfies
	$\lim_{n\to\infty}h_n=0$ and $\lim_{n\to\infty}(\log n)/(nh_n^q)=0$.

\item[(K4)] \textbf{Design points}: The covariates $\bX_1,\dots,\bX_n$ are iid random variables defined on a probability space $(\Omega,\mathcal{F},\mathbb{P})$ and valued on $[0,1]^q$ with  cdf $F_\bX$ and density $f_\bX$ that is bounded and bounded away from zero. Besides, they are independent of the errors $\varepsilon_1, \ldots, \varepsilon_n$.

\end{enumerate}

\begin{rem}
While consistency of estimators  defined in \eqref{estimatorIntegral}-\eqref{eq:integ} may be attained under mild assumptions of the boundary curve such as continuity of $g$ (see Neumeyer et al.\ (2019) for the proof in the univariate case and local polynomial estimator of order $0$), rates of convergence need stronger regularity assumption on $g$. Assumption (K1) is typical in such a framework and essentially means that $g$ is $\beta$-times differentiable and moreover a Lipschitz condition holds for the $\beta$ derivative. Further assumptions are also needed on the error distribution such as regular variation (K2), meaning that the distribution of the errors has polynomial tails. Assumptions (K1) and (K2) are common in the context of boundary models and are met in several papers, see  for instance Meister and  Rei{\ss} (2013), Jirak et al.\ (2014), M\"uller and Wefelmeyer (2010), Drees et al.\ (2019), Girard et al.\ (2013), H\"ardle et al.\ (1995) and Hall and Van Keilegom (2009); see also the book of de Haan and Ferreira (2006) for the applications and the motivation of heavy-tailed errors.
\end{rem}

\subsection{The fixed design case}
Similarly to the random design case, for $\bx \in [0,1]^q$, we consider the regression function estimator $\hat{g}$ defined as  
\begin{align}\label{estimatorIntegralFix}
\hat{g}_n(\bx):=\hat{g}(\bx):=p(\bx)
\end{align}
where $p$ is a polynomial of total degree $\beta^* \in \mathbb{N}_0$ and minimizes the local integral 
\begin{align}\label{eq:integFix}
\int_{[0,1]^q} p(\bt)I\{\|\bt-\bx\|\leq h_n\}d\bt
\end{align}
under the constraints $p(\bx_{i,n}) \geq Y_{i,n}$ for $\|\bx_{i,n}-\bx\|\leq h_n$. Here, $h_n > 0$ satisfies assumption (K3') below.\newline

We work for the fixed design case under (K1)-(K2) and the following modified assumptions (K3') and (K4'). 

\begin{enumerate}
	\item[(K3')] \textbf{Bandwidths}:  Let $(h_n)_{n \geq 0}$ be a sequence of positive bandwidths that satisfies
	$\lim_{n\to\infty}h_n=0$ and $\lim_{n \to \infty}(\log(n))/d_n=0$ with $d_n=d_n(1)$ from (K4').

\item[(K4')]  \textbf{Design points}: Let $I_n\subset [0,1]^q$  be a $q$-dimensional interval which is the cartesian product of one-dimensional closed intervals of length $dh_n$ with $d>0$. We assume that at least $d_n(d)$ design points lie in all such $I_n$.

\end{enumerate}

\begin{remark}
Assumptions (K3) and (K3') are common when analyzing the asymptotic behavior of such estimators of the form \eqref{estimatorIntegral} or  \eqref{estimatorIntegralFix}. In the equidistant fixed design framework we have $d_n(d)=nh_n^qd^q$ and (K3') equals (K3). Then for the univariate case $q=1$, we have $d_n=nh_nd$ and assumption (K3') turns to be assumption (H1) in Drees et al.\ (2019); see also assumption (A4') in Neumeyer et al.\ (2019).
\end{remark}

\begin{remark}
It is clear for assumption (K2) that the errors do not depend on the covariates in both cases of random and fixed design setups. Still, it has to be noticed that when dealing with the triangular scheme defined in  \eqref{model0}, the errors 
depend on $n$ too. This justifies the addition of the second index in $\varepsilon_{i,n}$. Indeed, the $i^{th}$ design point $\bx_{i,n}$ may vary with the sample size; see the construction of the fixed multivariate equidistant design in Section \ref{model-fixedequi} to be convinced. 
\end{remark}

\begin{remark}
The size of $\alpha$ in assumption (K2) is an important factor for the performance of the estimator defined in \eqref{estimatorIntegral}-\eqref{eq:integ} and \eqref{estimatorIntegralFix}-\eqref{eq:integFix} respectively. Simply speaking the smaller $\alpha>0$, the better the estimator. For $\alpha<2$ the error distribution is irregular and in this case the rate of convergence for the estimator is faster than the typical nonparametric rate, see the considerations below Theorem \ref{th:rateunif}.
In Figure \ref{figAlpha} we show some examples for different error distributions to highlight the effect of the size of $\alpha$. To simplify the presentation we restrict the display to $q=1$. 
\begin{figure}
\hspace{1.5cm}$\alpha=0.5$\hspace{3cm}$\alpha=1$\hspace{3cm}$\alpha=2$\hspace{3cm}$\alpha=3$\\
\begin{minipage}[c]{0.24\textwidth}
\includegraphics[width=\textwidth]{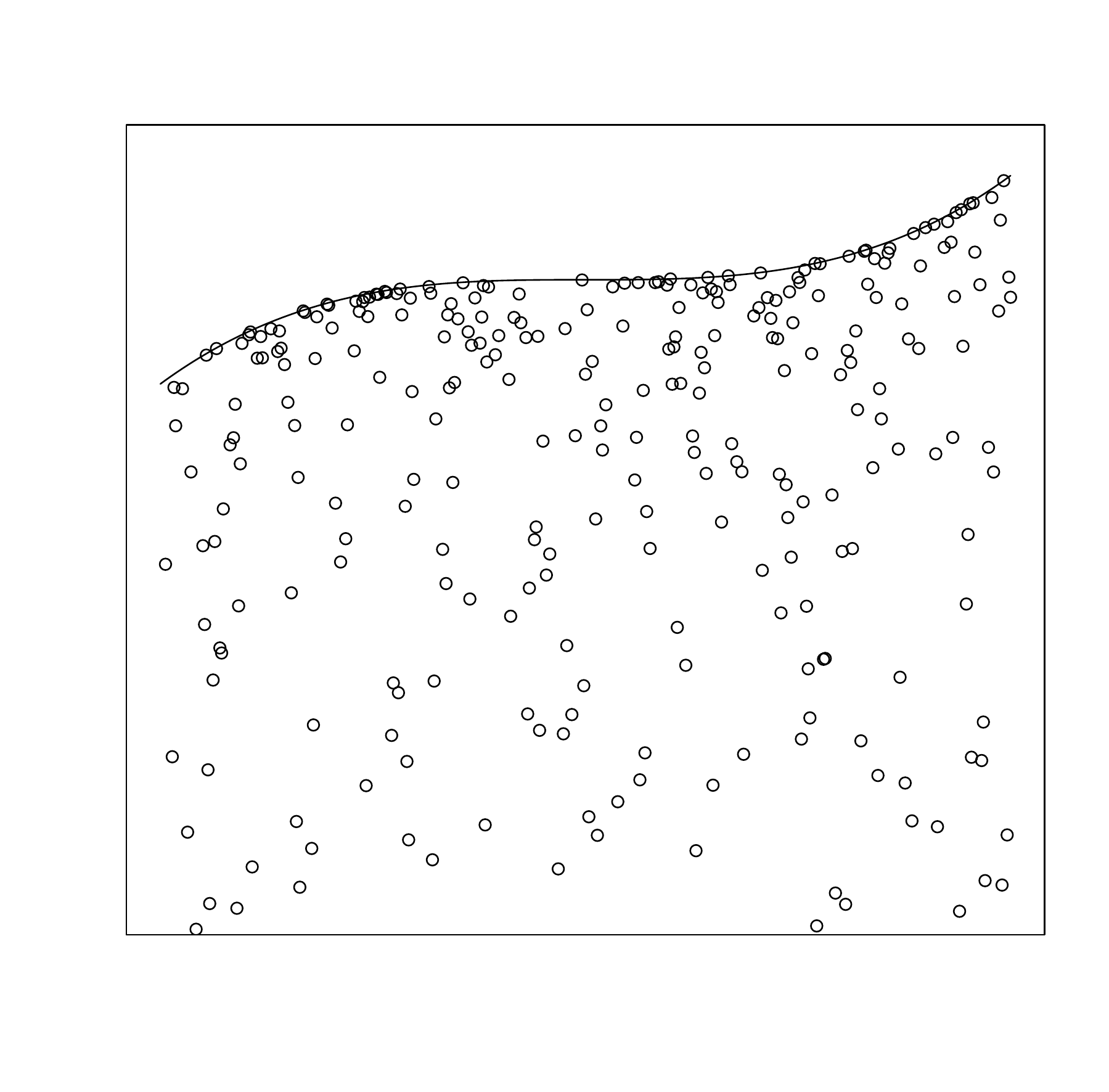}
\end{minipage}
\begin{minipage}[c]{0.24\textwidth}
\includegraphics[width=\textwidth]{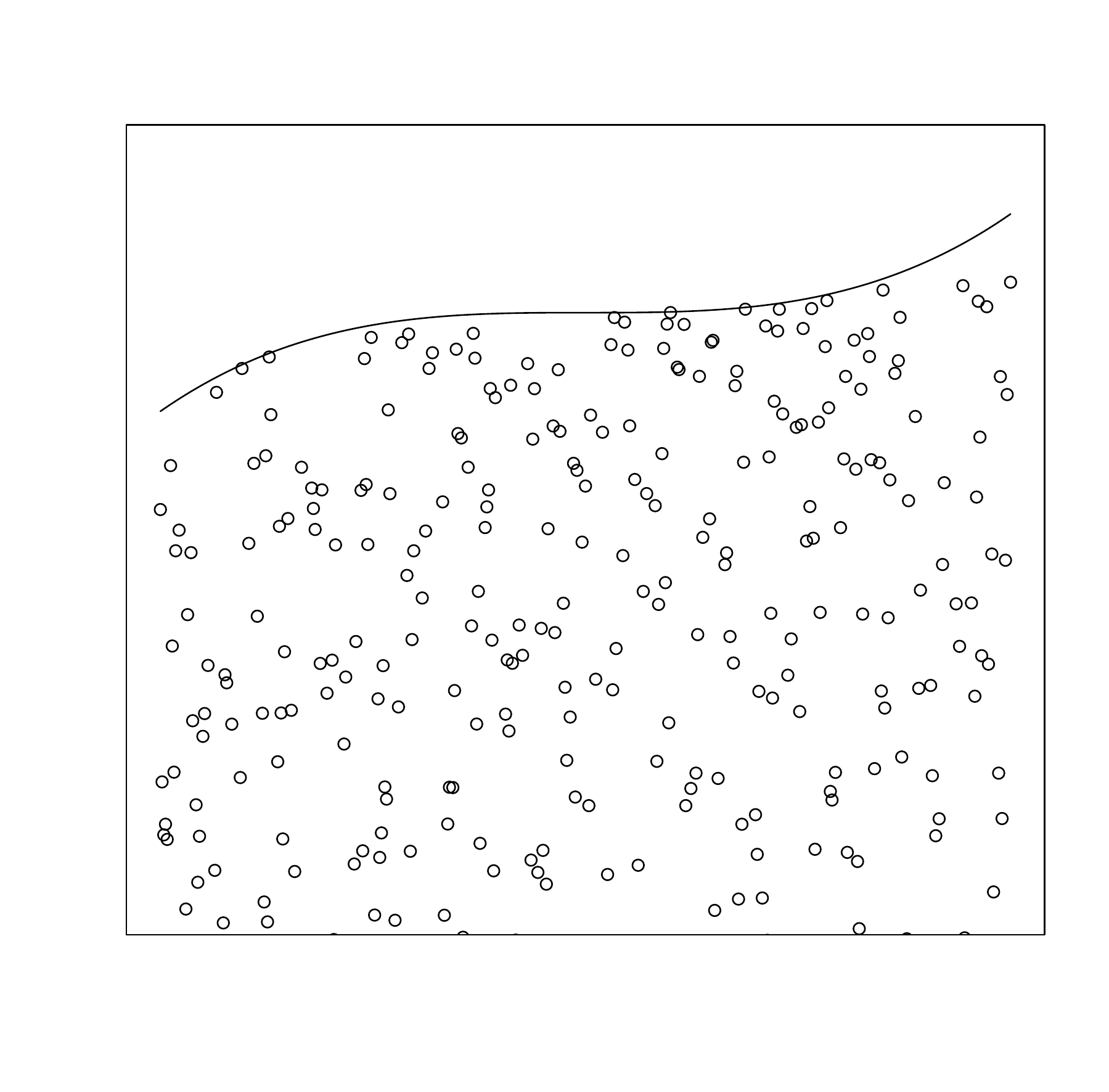}
\end{minipage}
\begin{minipage}[c]{0.24\textwidth}
\includegraphics[width=\textwidth]{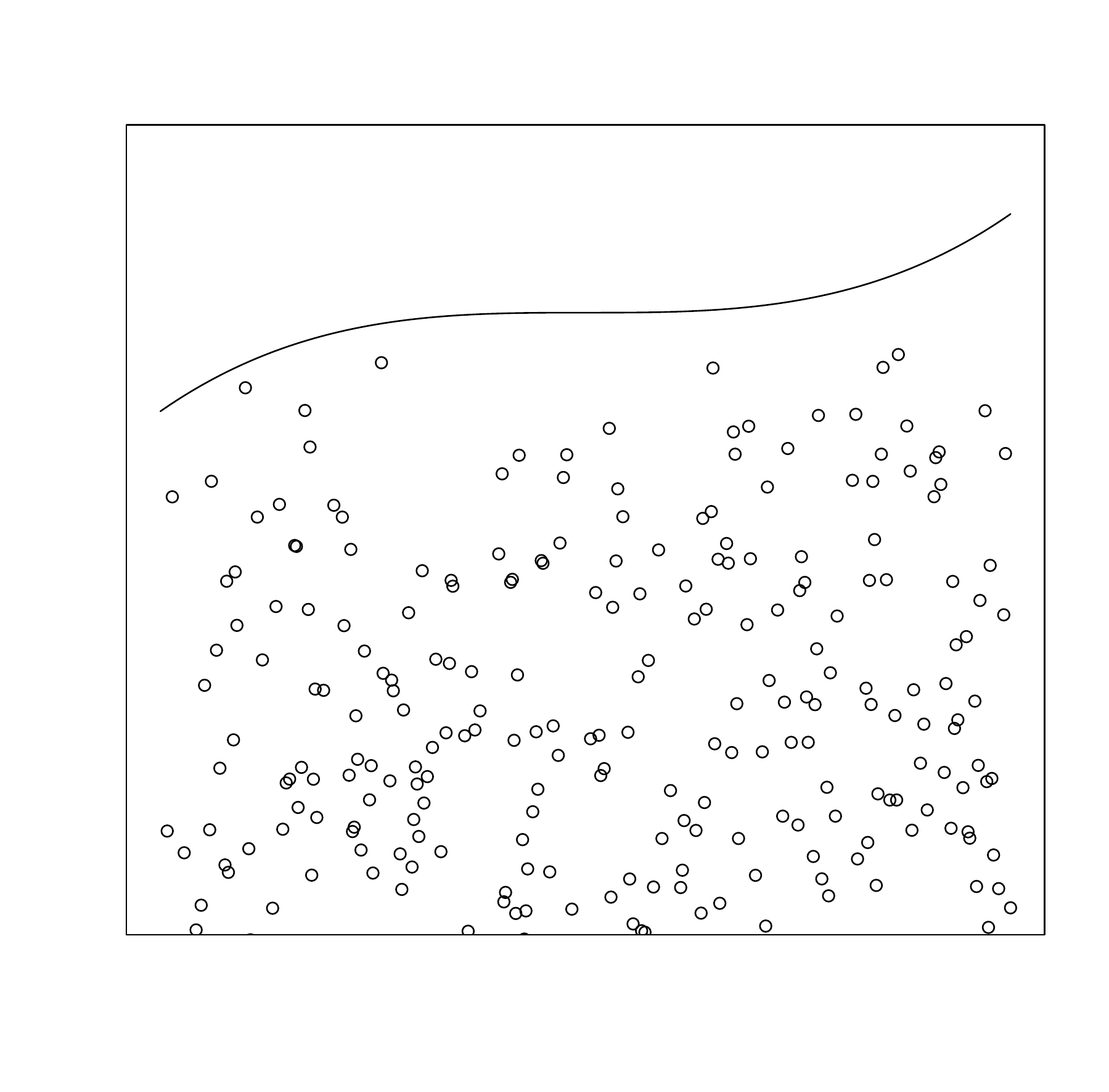}
\end{minipage}
\begin{minipage}[c]{0.24\textwidth}
\includegraphics[width=\textwidth]{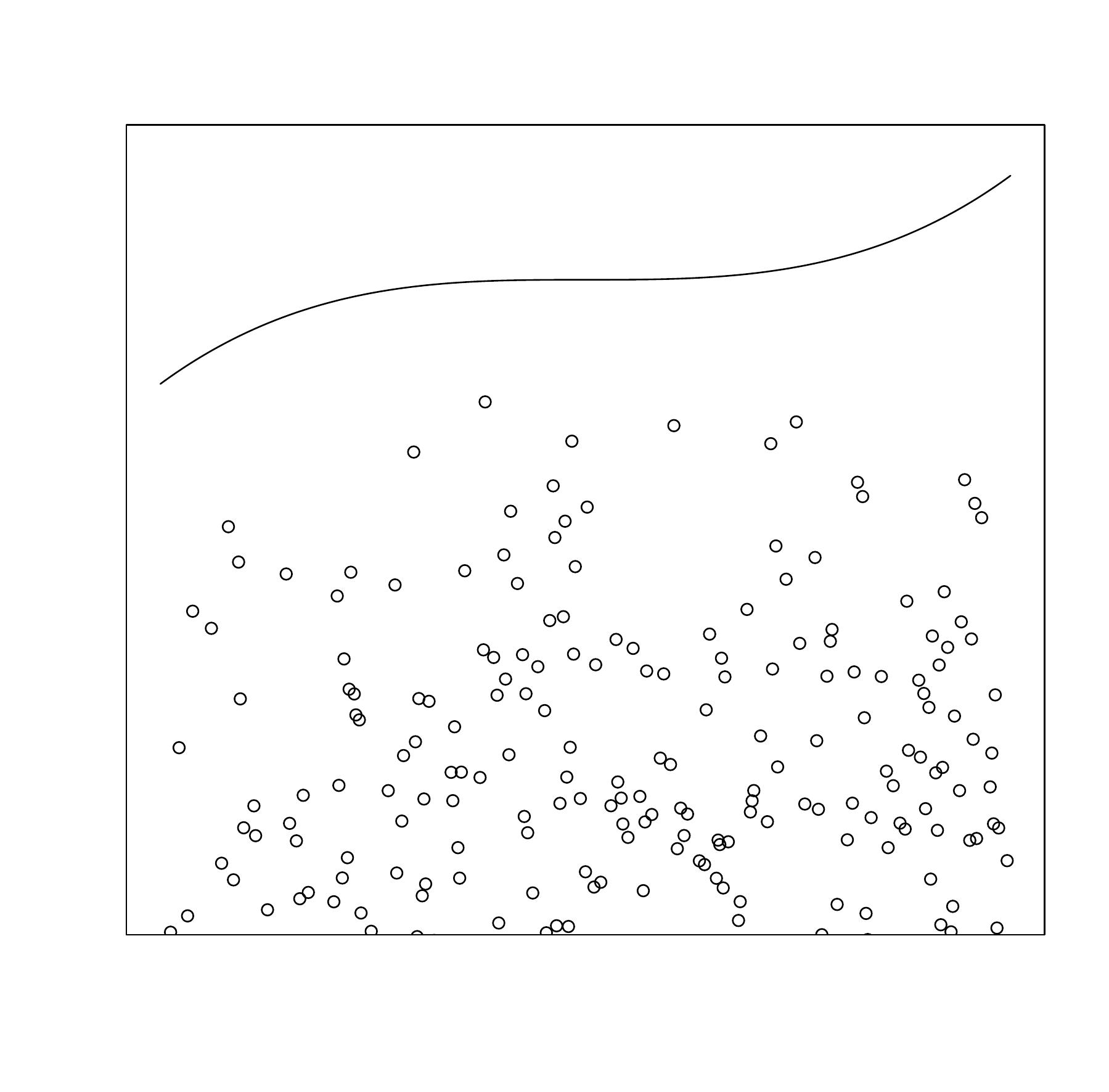}
\end{minipage}
\caption{Scatter plots of $(\frac in, Y_{i,n})$, $i=1,\ldots,n$ (fixed equidistant design for $q=1$) and the true regression function $g(\bx)=(\bx-0.5)^3+2$. The errors are Weibull distributed, s.\,t.\ $F(y)=\exp(-|y|^\alpha)I_{(-\infty,0)}(y)+I_{[0,\infty)}(y)$ for different values of $\alpha$.}\label{figAlpha}
\end{figure}
\end{remark}

\section{Main results}\label{main}

In this section we give the uniform consistency as well as the convergence rates of our estimator $\hat g$ separated for the two considered cases.

\subsection{The random design case}

In the next theorem, we provide the uniform consistency as well as the rate of convergence of the estimator of the regression function defined in \eqref{estimatorIntegral}-\eqref{eq:integ} for the random design case.

	\begin{theorem}\label{th:rateunif}
	Assume model \eqref{model-random} holds. If (K1) holds with $\beta  \in (0, \beta^*+1]$ and (K2)-(K4) are satisfied, then the estimator of the regression function defined in \eqref{estimatorIntegral}-\eqref{eq:integ} is uniformly consistent on  $[0,1]^q$ and we have 
	\begin{align*}
	\sup_{\bx \in [0,1]^q }|\hat{g}(\bx)-g(\bx)|=O(h_n^\beta)+O_\mathbb{P} \left(  \left(\frac{ \log(n)}{nh_n^q}\right)^{1/\alpha}\right).
	\end{align*}
\end{theorem}

\noindent\textbf{Proof:}
	See section \ref{appendix}. \newline

Note that the deterministic part $O(h_n^\beta)$ stems from the approximation of the regression function by a polynomial whereas the random part $O_\mathbb{P} \left(  \left(\frac{ \log(n)}{nh_n^q}\right)^{1/\alpha}\right)$ results from the observational error. Balancing the two error rates by setting $h_n = (\log(n)/n)^{\frac 1{\alpha\beta+q}}$ gives $\sup_{\bx \in [0,1]^q }|\hat{g}(\bx)-g(\bx)|=O_\mathbb{P}((\log(n)/n)^{\frac\beta{\alpha\beta+q}})$. For the case of an irregular error distribution, i.\,e.\ $\alpha\in(0,2)$, this rate improves upon the typical optimal rate $O_\mathbb{P}((\log(n)/n)^{\frac\beta{2\beta+q}})$ for the nonparametric estimation of mean regression functions in models with regular errors; see also the discussion in the simulation part.
	
\begin{remark}\label{rem:extension}
Theorem \ref{th:rateunif} extends the result of Drees et al.\ (2019), Theorem 2.2, to the multivariate random setting. Even in the univariate deterministic setting our result (Theorem \ref{th2:rateunif} established below) is an extension of the aforementioned Theorem since the convergence rate holds on the whole unit interval $[0,1]$ whereas in Drees et al.\ (2019) the result is restricted to $[h_n,1-h_n]$. The proof of the error rate that stems from the observational error follows along similar lines as the proof in Drees et al.\ (2019) but major adaptions are needed to deal with the multivariate and random case. The proof of the deterministic error rate that is due to approximating the boundary curve $g$ by a polynomial is based on the proof of Theorem 3.1 in Jirak et al.\ (2014). In the univariate equidistant fixed case that is treated in Drees et al.\ (2019) this Theorem can be applied directly whereas in the multivariate, possibly random case that is treated in the paper at hand the proof has to be intensely modified. Indeed, the original proof in Jirak et al.\ (2014) fully relies on the fundamental theorem of algebra, which states that every polynomial equation in one variable with complex
coefficients has at least one complex solution. As far as we know there is no possible extension of such a result for higher dimension hence moving to multidimensional covariates requires completely different arguments. This difficulty was already mentioned in Jirak et al.\ (2014), Remark 1. See the proof of Proposition \ref{prop:adapJirak} and especially the proof of Lemma \ref{inf-versus-sup} for our modification. This also gives an alternative proof of Theorem 3.1 in Jirak et al.\ (2014) and extends it to the multivariate and possibly random case.
\end{remark}

\subsection{The fixed design case}
 We give in the next theorem the uniform consistency as well as the rate of convergence of the estimator of the regression function defined in \eqref{estimatorIntegralFix}-\eqref{eq:integFix} for deterministic design points.

\begin{theorem}\label{th2:rateunif}
	Assume model \eqref{model0} holds. If (K1) holds with $\beta  \in (0, \beta^*+1]$ and (K2), (K3') and (K4') are satisfied, then the estimator of the regression function defined in \eqref{estimatorIntegralFix}-\eqref{eq:integFix} is uniformly consistent on $[0,1]^q$ and we have  
	\begin{align*}
	\sup_{\bx \in [0,1]^q }|\hat{g}(\bx)-g(\bx)|=O(h_n^\beta)+O_\mathbb{P} \left(  \left(\frac{\log(n)}{d_n}\right)^{1/\alpha}\right).
	\end{align*}
\end{theorem}
\noindent\textbf{Proof:}
See section \ref{appendix}. \newline

 When $\beta^*=0$ the H\"older class defined in Assumption (K1) reduces to the so-called class of \textit{$\beta$-Hölder uniformly continuous functions} with $\beta\in (0,1]$. In this framework, the boundary curve $g$ may be estimated by a local constant approximation 
	\begin{equation}\label{g_localconst}
	\hat g(\bx)=\max\{Y_{i,n}|i=1,\ldots,n \text{ with } \|\bx_{i,n}-\bx\|\leq h_n\}.
	\end{equation}
 This local constant approximation based estimator has been studied in Neumeyer et al.\ (2019) in the univariate setup for both random and fixed design points. Under the weaker assumption of continuity of the boundary curve $g$, they showed the uniform consistency of the estimator defined in \eqref{g_localconst} on the whole unit interval $[0,1]$. 
 Strengthening with the \textit{$\beta$-Hölder uniformly continuity} assumption and iid regularly varying innovations, we obtain from Theorem \ref{th2:rateunif} for the multivariate case, the uniform consistency as well as the rate of convergence for the deterministic design case on the whole unit interval $[0,1]^q$. We sum up in the following corollary.

\begin{cor}\label{prop:1}
	Assume (K2), (K3') and (K4') hold for model \eqref{model0} with $\beta$-Hölder uniformly continuous  function $g$. Then, the local constant approximation of the boundary curve $g$ defined in \eqref{g_localconst} is uniformly consistent on  $[0,1]^q$ and  we have  
	
	\begin{align*}	
	\sup_{\bx \in [0,1]^q} |\hat g(\bx)-g(\bx)| =O(h_n^\beta)+O_\mathbb{P} \left(  \left(\frac{ \log(n)}{d_n}\right)^{1/\alpha}\right).
	\end{align*}
\end{cor}
\begin{remark} Corollary \ref{prop:1} extends Remark 2.5 in Drees et al.\ (2019) to the multivariate setup and for non necessarily equidistant design points.
    
\end{remark}


\section{Simulations}\label{simulations}
To study the small sample behavior, we generate data according to the model
\[Y_i=0.5\cdot\sin\big(2\pi(X_i^{(1)}+\ldots+X_i^{(q)})\big)\ +\ 4(X_i^{(1)}+\ldots+X_i^{(q)})\ -\ \e_i\]
with $\bX_1,\ldots,\bX_n$ iid $\sim$ Unif$([0,1]^q)$ and $\e_1,\ldots,\e_n$ iid $\sim$ Exp$(1)$ for $q=2,3$ and sample sizes $n=10^q, 20^q, 30^q$. Thus in the simulated model the error distribution fulfills $\alpha=1$ and we consider different values of $\beta^*=0,1,2,3$ which means that we investigate local constant, local linear, local quadratic and local cubic approximation.

Since our estimator is computed by minimizing an integral over a polynomial the estimation procedure consists of solving a linear programming problem (compare to Remark \ref{Rem-solvePoly}).
The bandwidth is chosen as $h_n=n^{-\frac 1{\beta^*+1+q}}$ which corresponds (up to a log term) to the theoretically optimal bandwidth for $\alpha=1$ that balances the two error rates, see the considerations below Theorem \ref{th:rateunif}.

We run the same simulations with fixed equidistant design points $\bx_{1,n},\ldots,\bx_{n,n}$ as described in section \ref{model-fixedequi}. The results are very similar and thus we only present the results for the random design case. The discussion below also applies for the fixed design case.

 In Table \ref{table0.5} the results for 1000 replications are shown where we display the estimated mean squared error of our estimator $\hat g(0.5,0.5)$ ($\hat g(0.5,0.5,0.5)$ respectively). Except for the local constant approximation ($\beta^*=0$), the estimators perform well and become better as $\beta^*$ grows as well as $n$ grows; see the discussion below.

\begin{figure}
\hspace{1.5cm}$\beta^*=0$\hspace{3cm}$\beta^*=1$\hspace{3cm}$\beta^*=2$\hspace{3cm}$\beta^*=3$\\
\begin{minipage}[c]{0.24\textwidth}
\includegraphics[width=\textwidth]{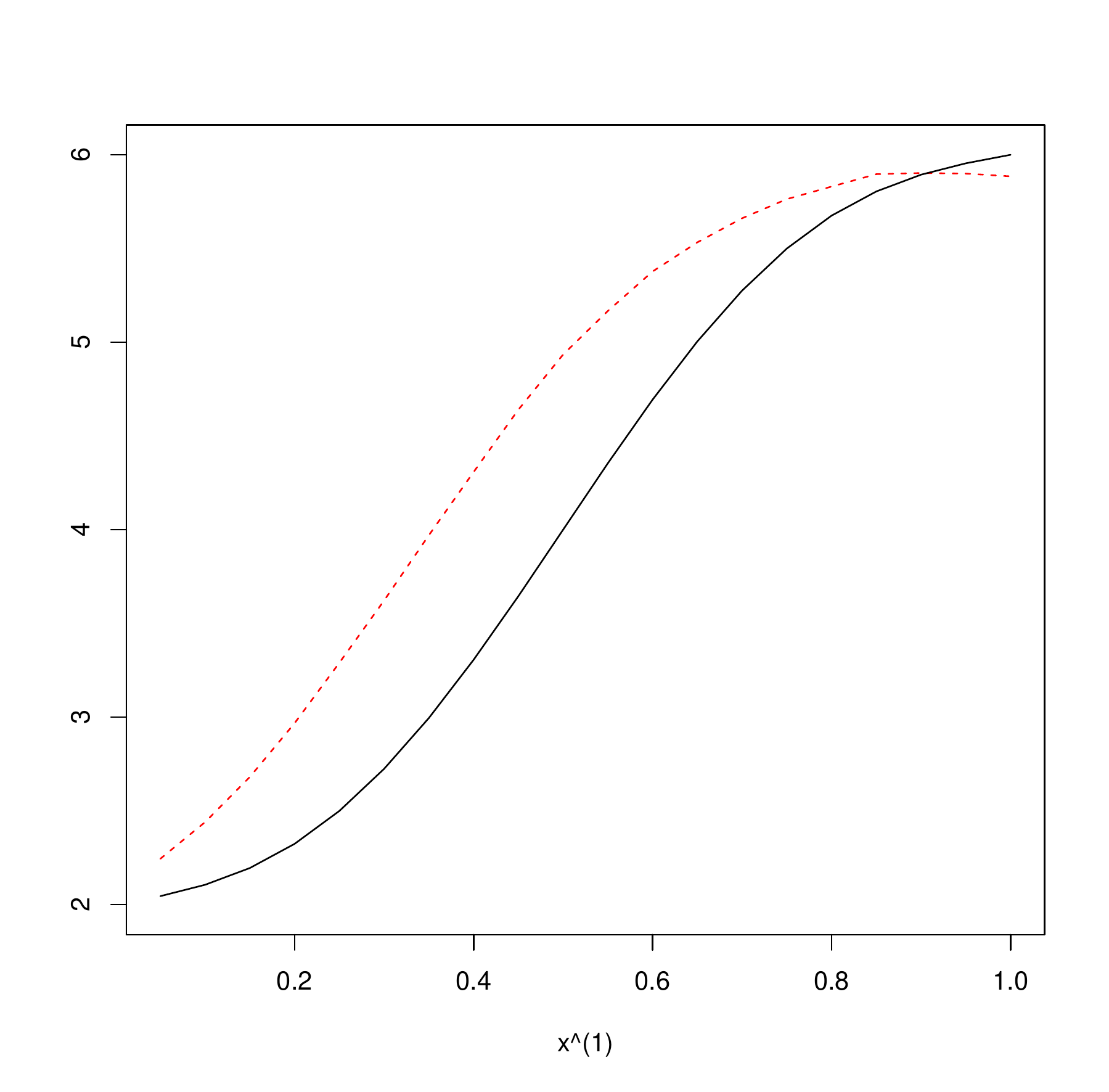}
\end{minipage}
\begin{minipage}[c]{0.24\textwidth}
\includegraphics[width=\textwidth]{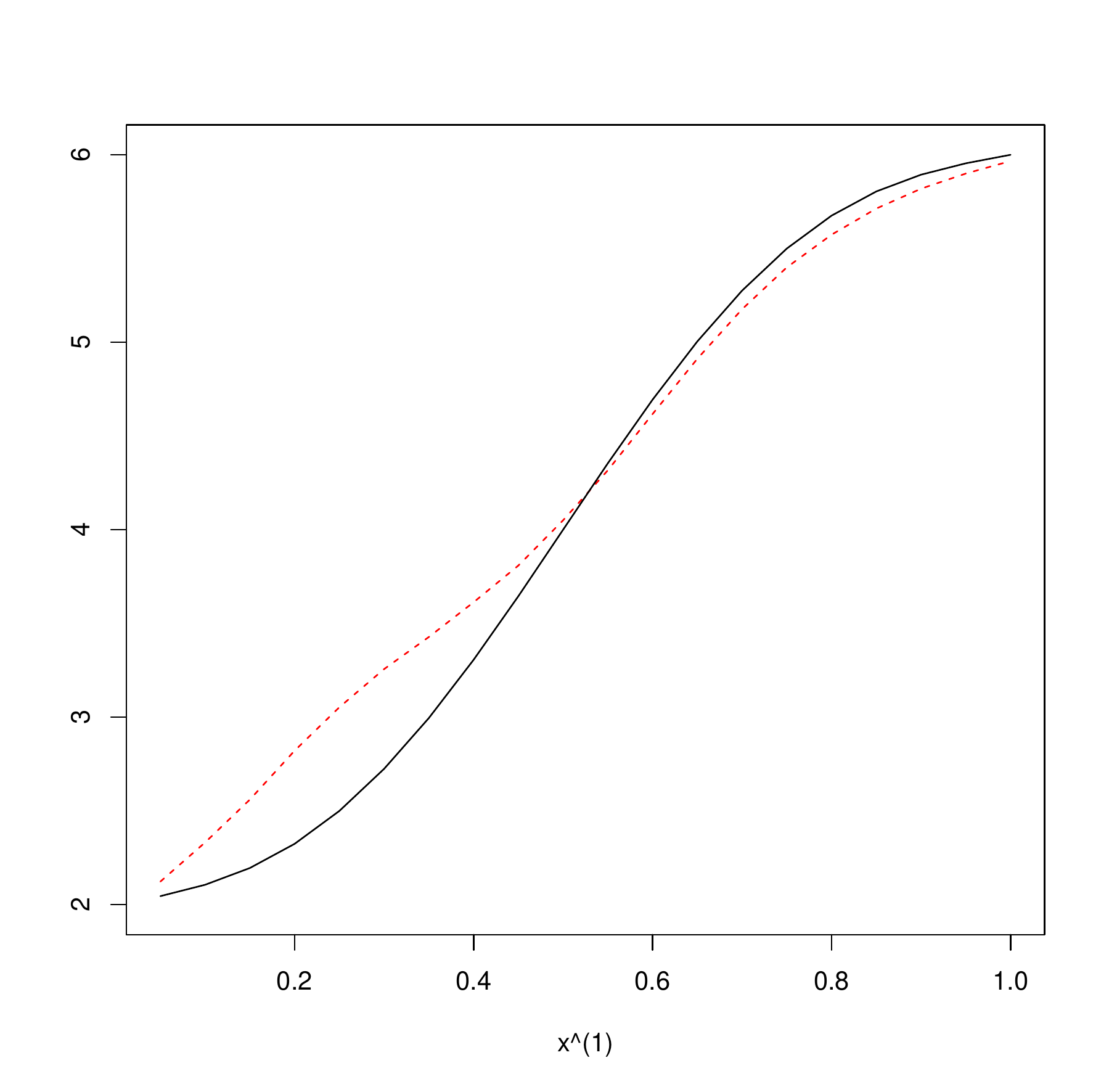}
\end{minipage}
\begin{minipage}[c]{0.24\textwidth}
\includegraphics[width=\textwidth]{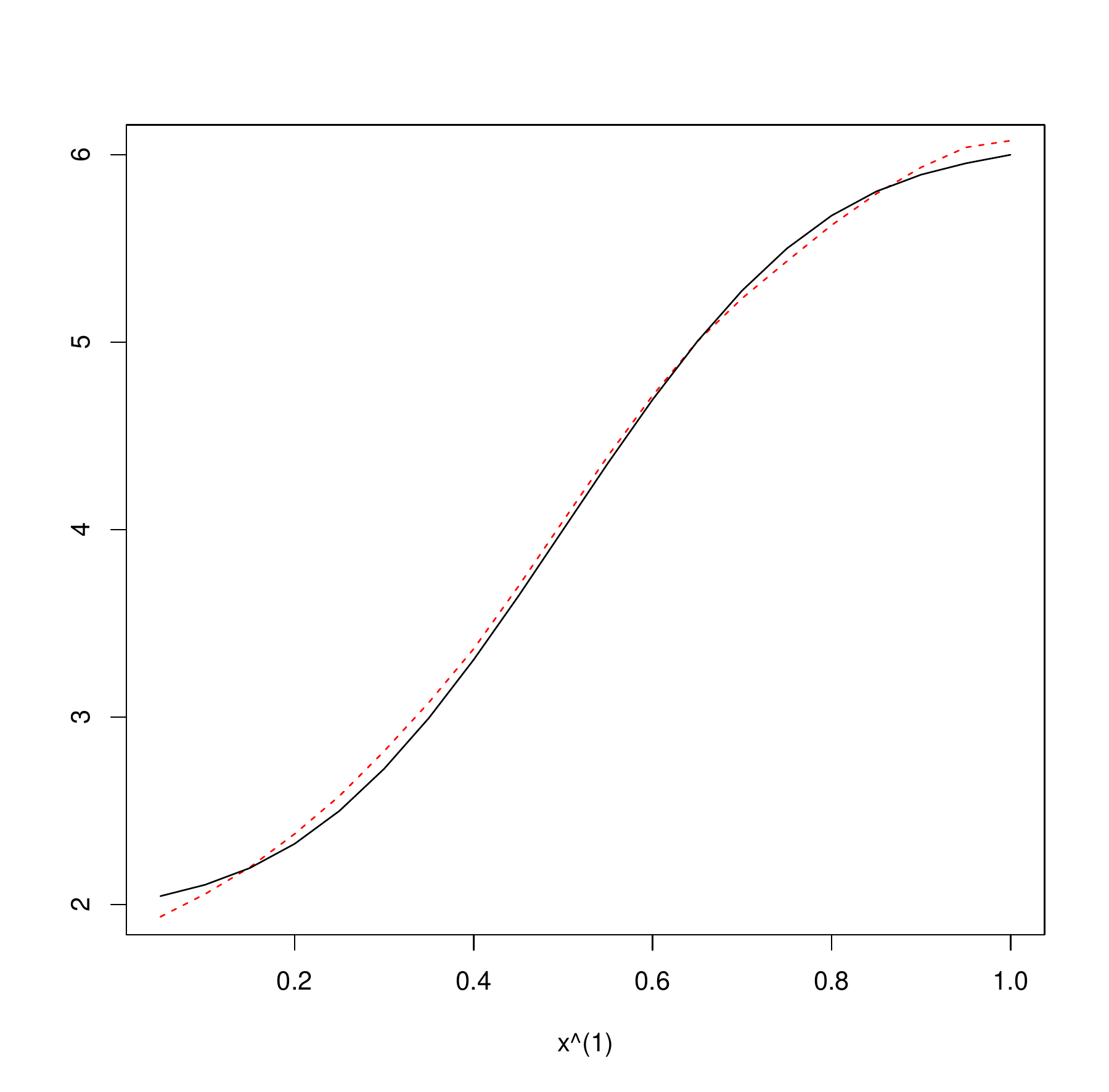}
\end{minipage}
\begin{minipage}[c]{0.24\textwidth}
\includegraphics[width=\textwidth]{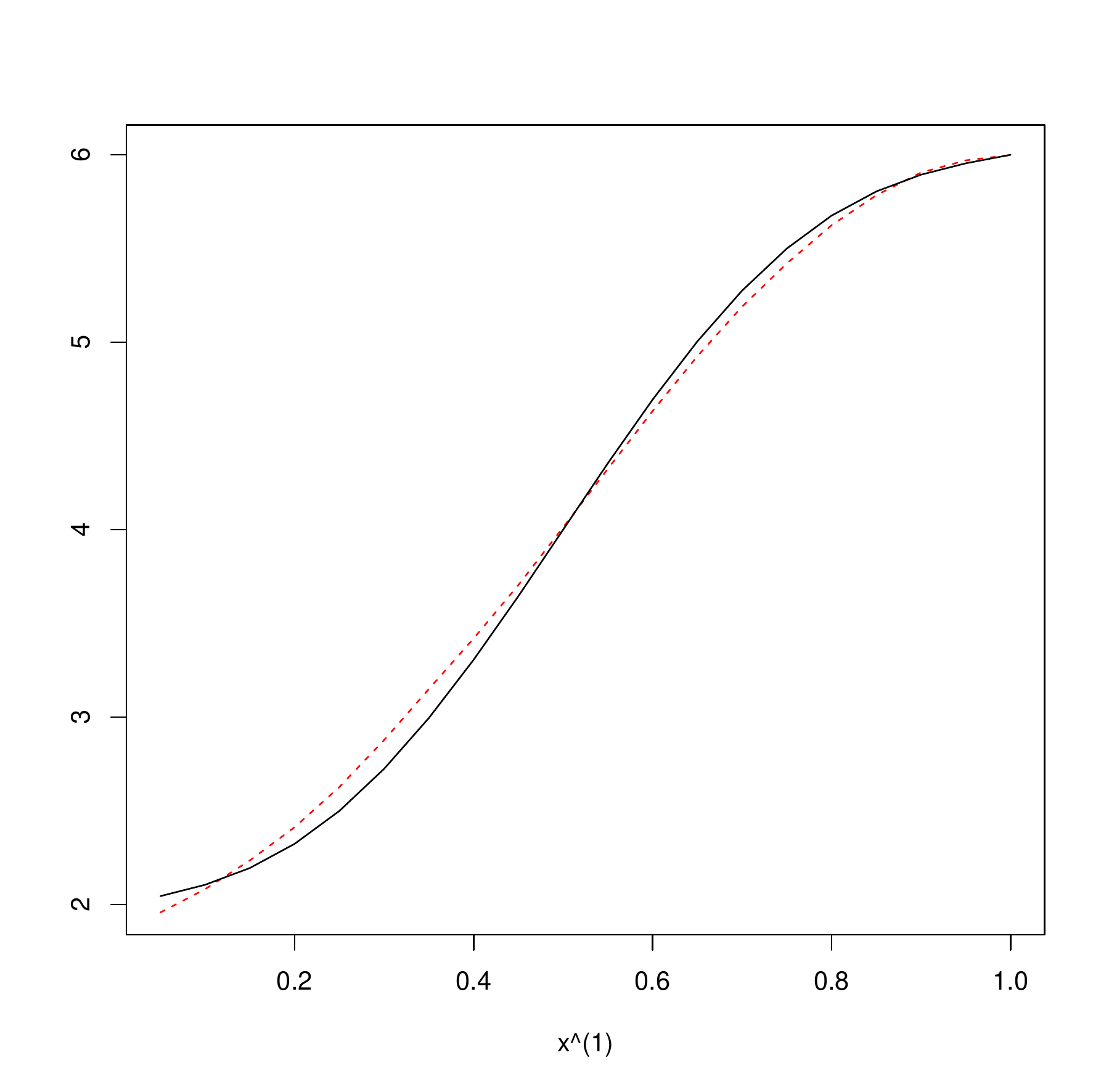}
\end{minipage}\\
\begin{minipage}[c]{0.24\textwidth}
\includegraphics[width=\textwidth]{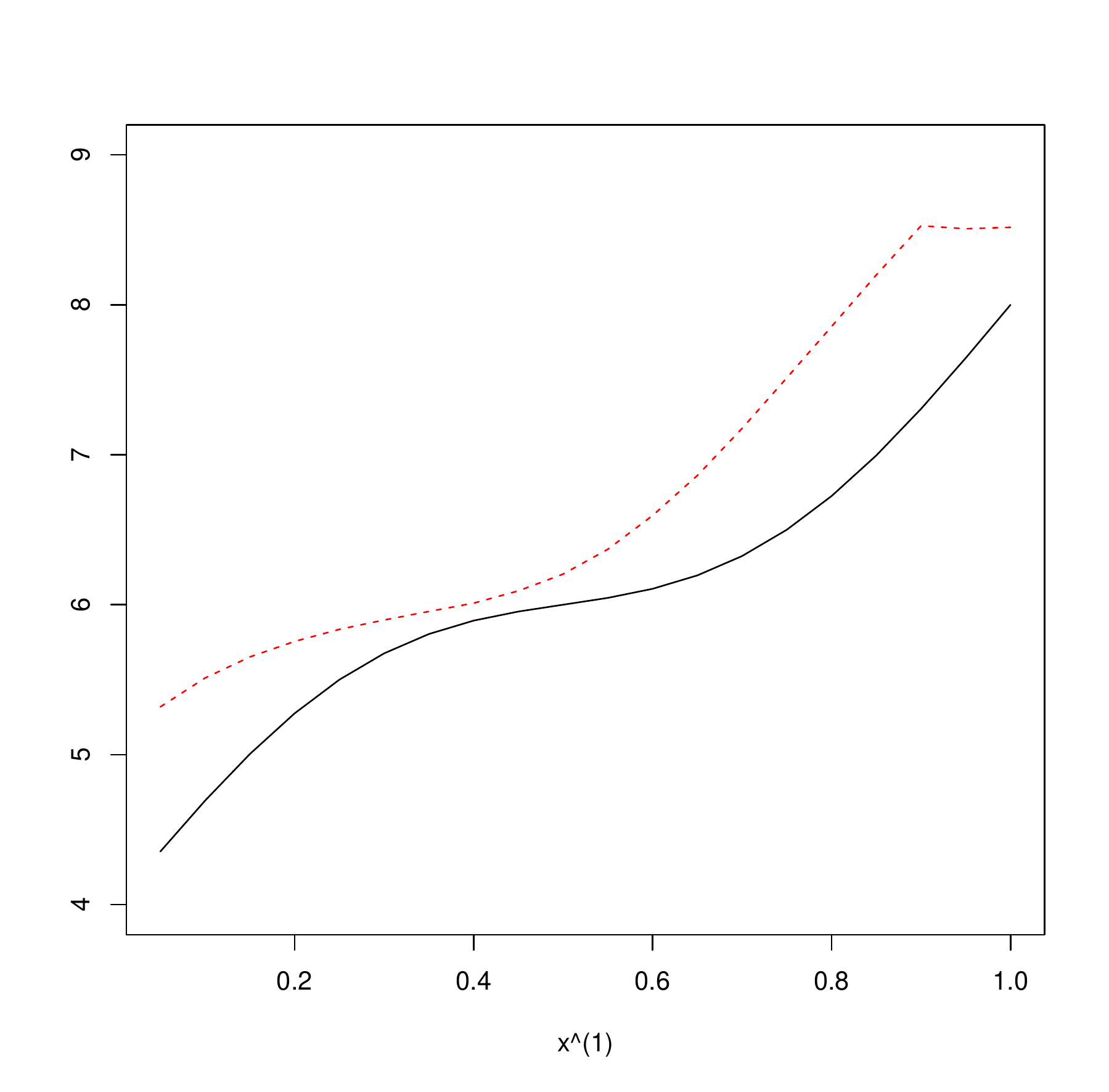}
\end{minipage}
\begin{minipage}[c]{0.24\textwidth}
\includegraphics[width=\textwidth]{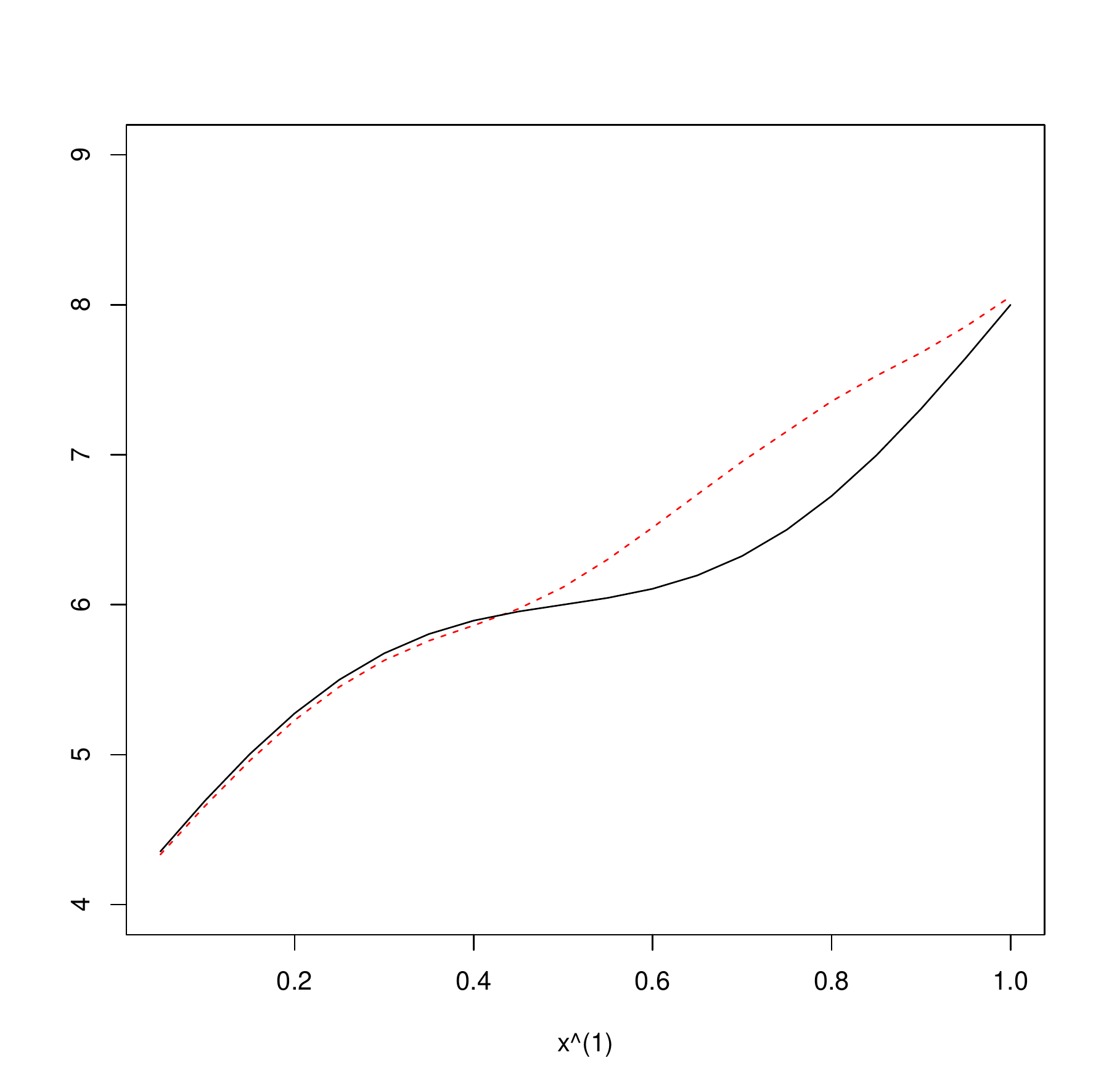}
\end{minipage}
\begin{minipage}[c]{0.24\textwidth}
\includegraphics[width=\textwidth]{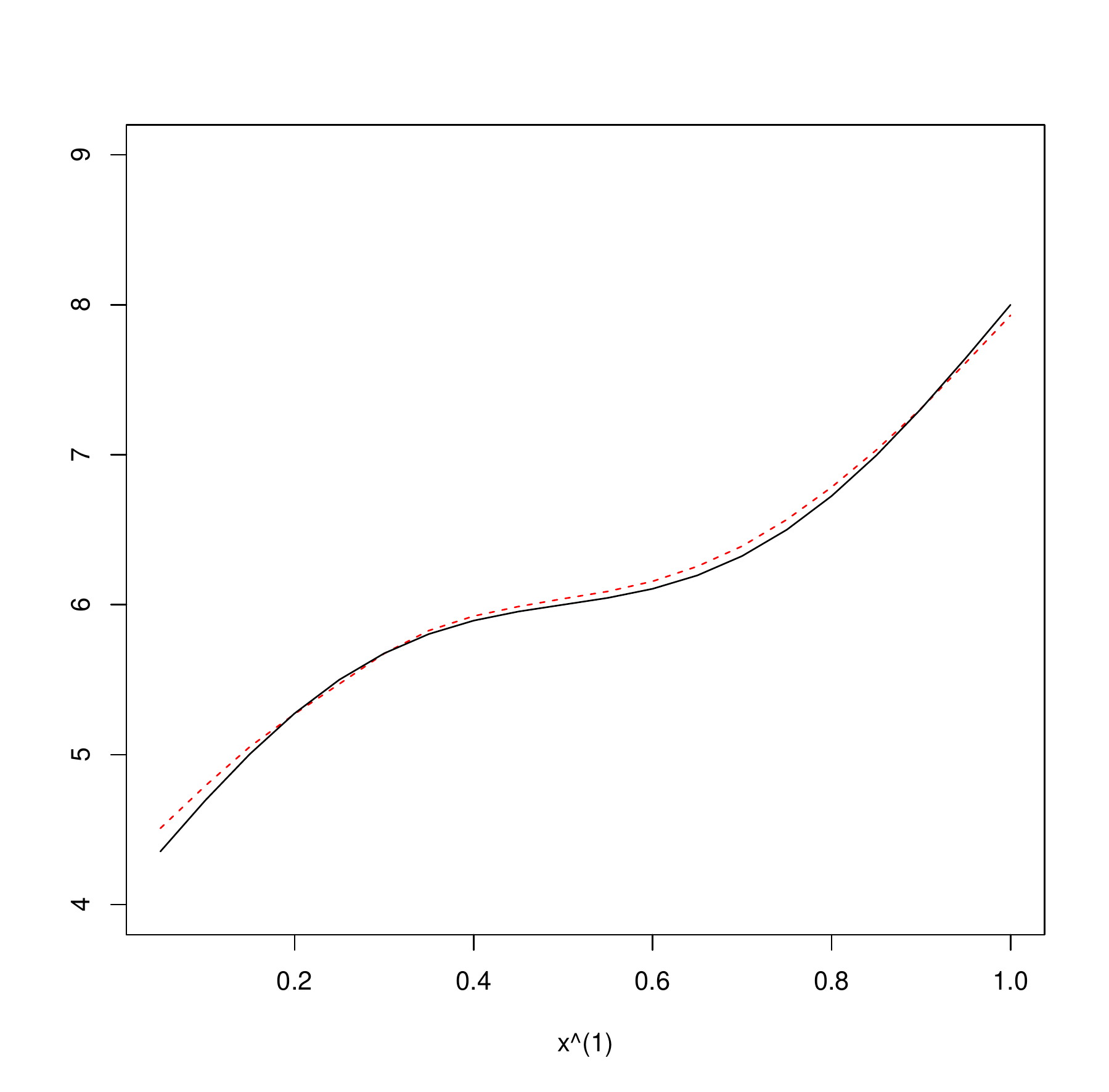}
\end{minipage}
\begin{minipage}[c]{0.24\textwidth}
\includegraphics[width=\textwidth]{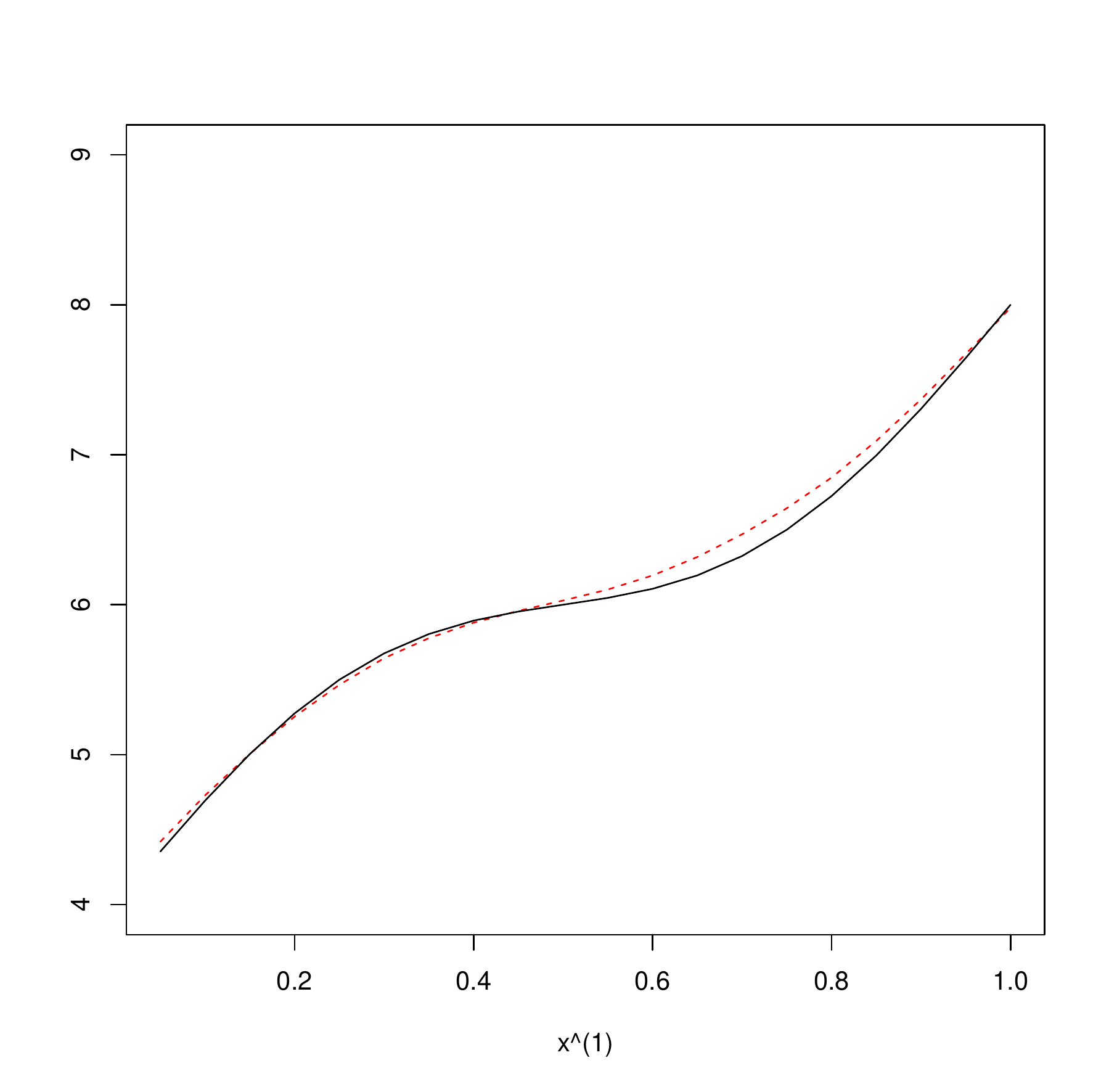}
\end{minipage}
\caption{True boundary curve (black solid line) and estimator $\hat g$ (red dashed line) for different values of $\beta^*$. The plot shows a cut trough the two dimensional surface (top) respectively three dimensional surface (bottom) of the function.}\label{figCurve}
\end{figure}

 In Figure \ref{figCurve} we show the true boundary curve and the estimated curve $\hat g$ in comparison. Since in dimension one the presentation of two curves in one plot is clearer than in higher dimensions we display a cut through the two respectively three dimensional surface of the functions. To be precise we plot the functions $x\mapsto 0.5\sin(2\pi(x+0.5))+4(x+0.5)$ and $x\mapsto \hat g(x,0.5)$ ($x\mapsto 0.5\sin(2\pi(x+1))+4(x+1)$ and $x\mapsto \hat g(x,0.5,0.5)$ respectively).
It can be seen that the approximation gets better as $\beta^*$ grows both for $q=2$ and $q=3$. For $\beta^*=2,3$ the approximation is very good, also in both considered cases of two and three dimensional covariates.
Exemplarily we show the results for $n=20^q$ since the effect is very similar for the other cases.

\begin{figure}
\begin{minipage}[l]{0.45\textwidth}
\includegraphics[width=\textwidth]{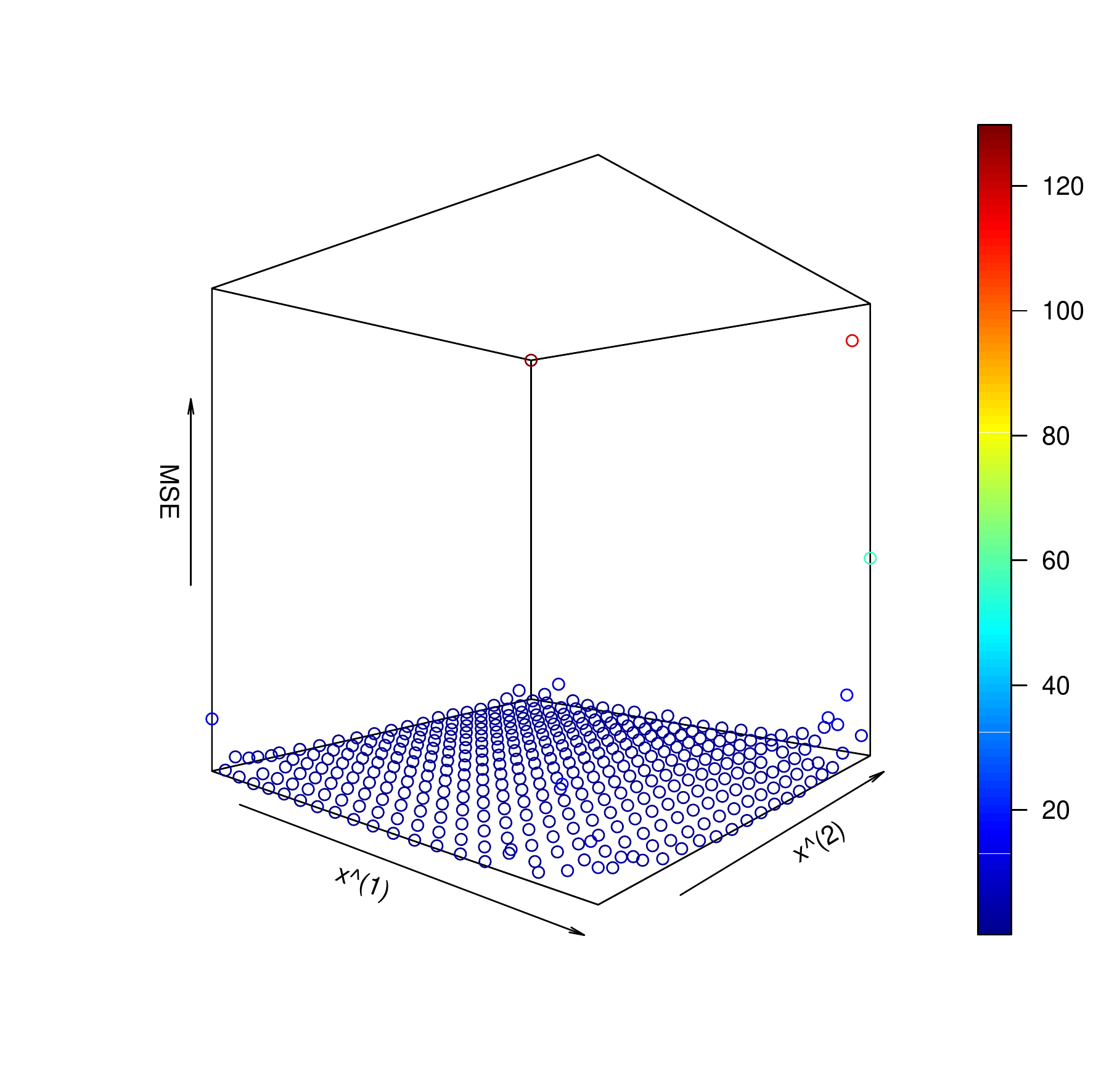}
\end{minipage}\qquad
\begin{minipage}[r]{0.45\textwidth}
\includegraphics[width=\textwidth]{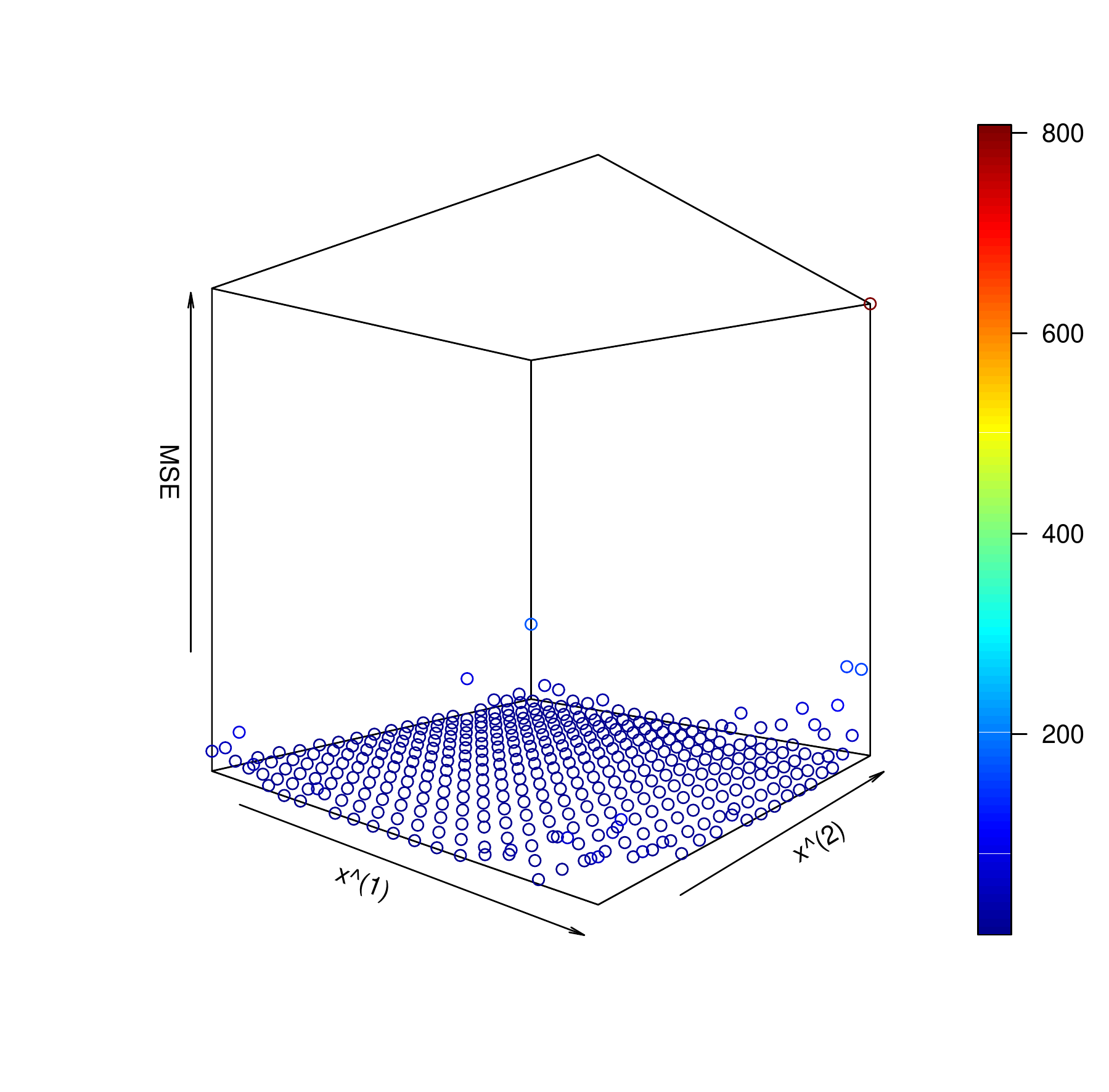}
\end{minipage}
\caption{Estimated mean squared error of $\hat g$ on $[0,1]^2$ with $\beta^*=2$ (left-hand side) and $\beta^*=3$ (right-hand side) and $n=10^2$.}\label{figMSE}
\end{figure}

 To evaluate the performance of the estimator on the whole interval $[0,1]^q$ we display in Table \ref{tableMean} the arithmetic mean of the estimated mean squared error of $\hat g(\bx_1),\ldots,\hat g(\bx_N)$ where $\bx_1,\ldots,\bx_N$ form a grid on $[0,1]^q$ with $N=20^q$.
Again, the performance is surprisingly poor for $n=10^2$ and get even worse for $\beta^*=2,3$; see the discussion below. In Figure \ref{figMSE} we show plots of the estimated mean squared error of $\hat g$ on $[0,1]^2$ for these cases. From the picture it can be deduced that the problem lies on the boundaries. The reason could be the smaller number of observations in this area which are of higher significance for larger $\beta^*$. 

While the performance of the estimators (except the two cases discussed before) are promising, we  highlight the following points:
\begin{itemize}

\item The higher the degree of the polynomial the more data is needed to attain satisfying estimation. 

\item For a point that lies in the center of $[0,1]^q$ the interval containing the observations for the estimation is of size $(2h_n)^q$, but if the point lies close to the boundary or at the boundary the interval gets smaller. Combined with the previous point, one has to be very careful when dealing with points at the boundary in small samples with polynomials of high degree.

\item As $n$ grows, the number of points within each interval grows.

\item The results improve as $\beta^*$ grows provided $n$ is large enough. The results improve as well as $n$ grows which are both expected effects that correspond to the theoretical result in Theorem \ref{th:rateunif}.

\item The local constant approximation (case $\beta^*=0$) seems to get worse as the dimension increases, in the case where not enough data points are provided or where the bandwidth is too large. One reason for this is that for each interval, the constant estimator always only fits one point (the maximum) whereas already in the linear estimator $(\beta^* = 1)$, the number of points fitted grows with the dimension (namely, $q+1$). For instance, for the bi-dimensional case with small sample size $n=10^2$ and $\beta^*=3$ the bandwidth $h_n$ is of order $n^{-1/6}\approx 0.46$ so $h_n$ may be too large: for a point $\bx$ far away from the boundaries, the resulting estimation $\hat g(\bx)$ using  local constant approximation consists in taking the maximum of the observations $Y_i$'s over (almost) the whole interval. When $n$ grows to $n=20^2$ and $n=30^2$ then the results gets better since the bandwidths shrink to $h_n \approx 0.14$ and $h_n \approx 0.10$.

\end{itemize}

Once mentioning that all of the issues discussed before fade away for $n$ large enough it has to be highlighted that bandwidth selection is a crucial point one has to be very careful with during the estimation procedure. From a theoretical perspective, it is clear that  $h_n$ drives the speed with which the boundary estimators converge; see Theorem 4.1. But even more importantly, when dealing with a small sample size in real data applications, the bandwidth impacts considerably the behavior of the estimator.

For the one-dimensional case, Jirak et al.\ (2014) propose a data-driven procedure to choose the bandwidth in any given application. The authors construct an adaptive estimator of the boundary of the curve based on a nested version of Lepski's method. They also investigate its performance in terms of pointwise and $L_q$ risks.
We believe that the method itself can be adapted to the higher-dimensional setting, because essentially it depends on the errors and the response variables, which are one-dimensional quantities.
The full theoretical analysis of such an adaption is a comprehensive topic and we leave it to future research.
Nevertheless, we sketch our adaption of the method below.

We start with constants $s \in (0,1)$ and $\rho>1$. The choice of these is decided by the user of the method. To avoid too large estimation steps,  $\rho$ should be close to $1$. The constant $s$ should not be too close to $1$.
We set $K = \lfloor \log_\rho (n^{1-s}) \rfloor$ and define the geometric sequence of bandwidths $h_k = h_0\rho^k$ for $k = 0,\dotsc,K+1$, where $h_0 = n^{s-1}$. For each $k$,
we define the estimator $\hat{g}_k$ as in Equations \eqref{estimatorIntegral} and \eqref{eq:integ}. We choose an optimal bandwidth $h_{\hat{k}}$ by setting
\[
\hat{k}:=\inf \{k=0, \ldots,K : \exists l \leq k : \Vert \hat{g}_{k+1}-\hat{g}_l \Vert  > \hat{\zeta_l} + \hat{\zeta}_{k+1}\} \wedge K
\]
and choosing the corresponding $h_{\hat{k}}$.
Here, the thresholds $\hat\zeta_{k}$ can be estimated from the one-dimensional observations $(Y_i)_{i=1,\dotsc,n}$ as explained in Jirak et al.\ (2014) and sketched below.
Finally, the required data-driven estimator is defined as $\hat{g}_{\text{driven}}:=\hat{g}_{\hat{k}}$.


The estimation of the thresholds $\zeta_k$ requires first the estimation of the extremal indices involved in the error distribution; namely the tail index denoted $\alpha$ in the present paper as well as the coefficients involved in the second order conditions of regular variation, see Equations (1.3) and (1.4) and Assumption 3.1 in Jirak et al.\ (2014) for more details. Since these coefficients are not directly related to the dimension of the covariates, their estimation could be done using similar arguments as in Jirak et al.\ (2014) i.e., with a Hill-type estimator; see their Section 3.2 and more specifically Equation (3.10) for the precise definition of the estimator.
This part requires an investigation of the tail distribution of the error and involves some additional technical assumptions; see for instance Equation (3.17) and Assumption 10.1 in their paper. In sum, we believe that the estimation of the critical thresholds in the multivariate setting may be attained in the same fashion as in their Sections (3.3) and (3.4).

Note that different approaches to this problem may be considered such as a theoretically lighter cross validation approach or a bootstrap method, we refer to Hall and Park (2004) for more details.

\begin{table}
\begin{center}
\begin{tabular}{l||c|c|c}
 $q=2$ & $n=10^2$ & $n=20^2$ & $n=30^2$\\
 \hline
$\beta^*=0$ & $1.56$ & $0.95$ & $0.62$\\
$\beta^*=1$ & $0.02$ & $0.006$ & $0.002$\\
$\beta^*=2$ & $0.07$ & $0.01$ & $0.005$\\
$\beta^*=3$ & $0.02$ & $0.003$ & $0.001$
\end{tabular}
\quad 
\begin{tabular}{l||c|c|c}
 $q=3$ & $n=10^3$ & $n=20^3$ & $n=30^3$\\
 \hline
$\beta^*=0$ & $0.62$ & $0.06$ & $0.02$\\
$\beta^*=1$ & $0.06$ & $0.02$ & $0.005$\\
$\beta^*=2$ & $0.03$ & $0.003$ & $0.001$\\
$\beta^*=3$ & $0.02$ & $0.001$ & $0.0003$
\end{tabular}
\end{center}
\caption{Estimated mean squared error for $\hat g(0.5,0.5)$ (left-hand side) and for $\hat g (0.5,0.5,0.5)$ (right-hand side)}
\label{table0.5}
\end{table}

\begin{table}
\begin{center}
\begin{tabular}{l||c|c|c}
 $q=2$ & $n=10^2$ & $n=20^2$ & $n=30^2$\\
 \hline
$\beta^*=0$ & $0.79$ & $0.37$ & $0.23$\\
$\beta^*=1$ & $0.22$ & $0.06$ & $0.04$\\
$\beta^*=2$ & $1.4$ & $0.03$ & $0.009$\\
$\beta^*=3$ & $7.47$ & $0.05$ & $0.01$
\end{tabular}
\quad 
\begin{tabular}{l||c|c|c}
 $q=3$ & $n=10^3$ & $n=20^3$ & $n=30^3$\\
 \hline
$\beta^*=0$ & $1.18$ & $0.52$ & $0.3$\\
$\beta^*=1$ & $0.15$ & $0.08$ & $0.05$\\
$\beta^*=2$ & $0.09$ & $0.008$ & $0.002$\\
$\beta^*=3$ & $0.56$ & $0.008$ & $0.002$
\end{tabular}
\end{center}
\caption{Arithmetic mean of estimated mean squared errors for $\hat g$ on $[0,1]^2$ (left-hand side) and on $[0,1]^3$ (right-hand side)}
\label{tableMean}
\end{table}

\section{Proofs of the main results}\label{appendix}

\subsection{Proof of Theorem \ref{th:rateunif}}
The proof of Theorem \ref{th:rateunif} is a direct application of Propositions \ref{prop:adapJirak} and \ref{min-eps-rand} stated below and follows along similar lines as the proof of Theorem 2.2 in Drees et al.\ (2019).
\begin{proposition}\label{prop:adapJirak}
	Assume that model \eqref{model-random} holds under (K2) and (K4) and consider the regression function estimator $\hat{g}$ defined in \eqref{estimatorIntegral}-\eqref{eq:integ} where $g$ fulfills condition (K1) for some $\beta \in (0, \beta^* +1]$ and some $c_g \in [0,c^*]$. Then, there exist $C_{\beta^*,q,c^*},C_{\beta^*,q}$ and a natural number $J_{\beta^*, q}$, which depend only on the respective subscripts such that for all $\bx \in [0,1]^q $
	\begin{align*}
	|\hat{g}(\bx)-g(\bx)|\leq C_{\beta^*,q,c^*}h_n^\beta +  C_{\beta ^*,q}  \max_{\bj\in\{1,\ldots ,2 J_{\beta^*, q}\}^q}  \left(    
	\min_{\bX_i \in \bx+h_n I_{\bj}}	|\varepsilon_{i}|   \right),
	\end{align*}
	where $\bx+h_nI_\bj=\{(x^{(1)}+h_n I_{j_1})\cap[0,1]\}\times\ldots\times \{(x^{(q)}+h_nI_{j_q})\cap[0,1]\}$ with $I_k=[-1+(k-1)/J_{\beta^*, q}, -1 + k/J_{\beta^*, q}]$.
\end{proposition}

\noindent\textbf{Proof:}	We first highlight that throughout the proof the design points $\bX_i:(\Omega, \mathcal{F}, \mathbb{P}) \longrightarrow ([0,1]^q, \mathcal{B}[0,1]^q)$ are random elements. To make the reading easier, we omit the script $\omega$ but the proof has to be understood $\omega$-wise that is for any realisation $\bX_i(\omega), \omega \in \Omega$. Besides, unless it is specified otherwise, $n \in \mathbb{N}$ is arbitrary. \newline

	Let  $\bx \in [0,1]^q$ be fixed and for $n \geq 1$ set $I_n^*=\{[x^{(1)}-h_n,x^{(1)}+h_n] \cap [0,1]\}\times\ldots\times\{[x^{(q)}-h_n,x^{(q)}+h_n] \cap [0,1]\}$. The idea of proof is based on Theorem 3.1 in Jirak, Meister and Rei{\ss} (2014) but comprehensive adaptions are needed to deal with the multivariate case. We consider random design points satisfying assumption (K4) where the Riemann approximation in the aforementioned paper 
	 is replaced by the integral defined in \eqref{eq:integ}.\newline  
	 
\noindent This means that we consider the coefficients $(\hat b_\bj)_\bj$ for all multiindices $\bj$ with $|\bj|\in\{0,\ldots,\beta^*\}$ which minimize the objective function 
	\begin{align}\label{locInt}
	S(\bx,(b_\bj)_\bj)& = \int_{I_n^*}\sum_{\bj\in\en_0^q:|\bj|\leq\beta^*}b_\bj(\bt-\bx)^\bj d\bt
	\end{align}
	under the constraint $Y_i \leq \sum_{\bj\in\en_0^q:|\bj|\leq\beta^*}b_\bj(\bX_i-\bx)^\bj$ for all $i$ with $\bX_i\in I_n^*$.	
\newline

\noindent Now, a Taylor-Lagrange development up to the order $\langle \beta \rangle$ of $g$ around $\bx$ yields  
	\begin{align*}
	g(\bX_i)=\sum_{|\bj|\leq\langle \beta \rangle -1}\frac{D^\bj g(\bx)}{\bj!}(\bX_i-\bx)^\bj  +\sum_{|\bj|=\langle\beta\rangle}\frac{D^\bj g(\bx+ \theta(\bX_i-\bx))}{\bj!}(\bX_i-\bx)^\bj 
	\end{align*}
		with $\theta \in [0,1]$. Then we have 
\begin{align*}
g(\bX_i) & =\sum_{|\bj|\leq \langle \beta \rangle -1}\frac{D^\bj g(\bx)}{\bj!}(\bX_i-\bx)^\bj +\sum_{|\bj|=\langle\beta\rangle}\frac{D^\bj g(\bx)}{\bj!}(\bX_i-\bx)^\bj +\sum_{|\bj|=\langle\beta\rangle}\frac{D^\bj g(\bx+ \theta(\bX_i-\bx))}{\bj!}(\bX_i-\bx)^\bj \\
& - \sum_{|\bj|=\langle\beta\rangle}\frac{D^\bj g(\bx)}{\bj!}(\bX_i-\bx)^\bj \\
& = \sum_{|\bj|\leq \langle \beta \rangle }\frac{D^\bj g(\bx)}{\bj!}(\bX_i-\bx)^\bj + r_{\langle \beta \rangle}(\bX_i,\bx)
\end{align*}
where $r_{\langle \beta \rangle}$ is the remainder term defined for $\bt,\bx \in [0,1]^q$ by
	\begin{align*}
r_{\langle \beta \rangle}(\bt,\bx)=  \sum_{|\bj|=\langle\beta\rangle}\frac{D^\bj g(\bx+ \theta(\bt-\bx))-D^\bj g(\bx)}{\bj!}(\bt-\bx)^\bj  .
\end{align*}
By assumption (K1) we can make use of the Hölder property \eqref{def:holder} and get
\beq
r_{\langle \beta \rangle}(\bt,\bx)&\leq&c_g\|\bt-\bx\|^{\beta-\langle\beta\rangle} \sum_{|\bj|=\langle\beta\rangle}\frac{(\bt-\bx)^\bj}{\bj!}\\
&\leq&c_g\|\bt-\bx\|^{\beta-\langle\beta\rangle} \sum_{|\bj|=\langle\beta\rangle}\frac{\|\bt-\bx\|^{\langle\beta\rangle}}{\bj!}\\
&=&c_gc_1(\beta,q)\|\bt-\bx\|^\beta
\eeq
with some constants $c_g,c_1(\beta,q)<\infty$.\newline

\noindent Consider now that $b_\bj$ are the Taylor coefficients such that 
\begin{align*}
b_\bj=\left\{
\begin{array}{l}
D^\bj g(\bx)/\bj! \ \ \ \text{if} \ \ \  \ \ \ |\bj| \leq \langle \beta \rangle\\
0 \ \ \ \ \ \ \ \ \ \ \ \ \ \ \text{if} \ \ \  \ \ \ |\bj| > \langle \beta \rangle.
\end{array}
\right.
\end{align*}
With this we define the following two quantities 
	\begin{align*}
	\bigtriangledown_{i}=g(\bX_i)-\sum_{|\bj|\leq\beta^*}b_\bj(\bX_i-\bx)^\bj, \ \ \  \ \ \ i=1,\ldots,n
	\end{align*}
and
\begin{align*}
	\bigtriangledown_{n}^*:= \sup_{\bt\in I_n^*}|r_{\langle \beta \rangle}(\bt,\bx)|.
\end{align*}
	Then, one may rewrite the data points $Y_{i}$ in the model \eqref{model-random} as 
	\begin{align*}
	Y_{i}=\sum_{|\bj|\leq\beta^*}b_\bj(\bX_i-\bx)^\bj + \varepsilon_i +\bigtriangledown_{i}, \ \ \ i=1,\ldots,n,
	\end{align*}
and from what precedes, under the constraint $\bX_i\in I_n^*$, it follows that
	\begin{align}\label{infDelta}
	\bigtriangledown_{i}\leq \bigtriangledown_{n}^* \leq  c_g c_1(\beta,q) h_n^\beta.
	\end{align}
%
%
Since the errors $\varepsilon_{i},i=1, \ldots,n$ are non-positive, we have  for any $\bX_i \in  I_n^* $
	\begin{align*}
	Y_{i} & = \sum_{|\bj|\leq\beta^*}b_\bj(\bX_i-\bx)^\bj + \varepsilon_{i} +\bigtriangledown_{i}\\
	& \leq \sum_{|\bj|\leq\beta^*}b_\bj(\bX_i-\bx)^\bj  +\bigtriangledown_{i}\\
	& \leq \sum_{|\bj|\leq\beta^*}b_\bj(\bX_i-\bx)^\bj   +\bigtriangledown_{n}^*, \ \ \ i=1, \ldots,n.
	\end{align*}
Thus, since the coefficients $(\hat{b}_\bj)_\bj$  minimize the local integral \eqref{locInt}, we have 
	\begin{align}\label{ineq:1}
	\int_{I_n^*} \sum_{|\bj|\leq\beta^*}\hat{b}_\bj(\bt-\bx)^\bj d\bt & \leq \int_{I_n^*} \left( \sum_{|\bj|\leq\beta^*}b_\bj(\bt-\bx)^\bj   + \bigtriangledown_{n}^* \right) d\bt
	\end{align}
	for all $n$.
Define now the polynomial
	\begin{align*}
	Q(\bt)=\sum_{|\bj|\leq\beta^*}(b_\bj-\hat{b}_\bj)(\bt-\bx)^\bj+ \bigtriangledown_{n}^*, \ \ \ \bt \in I_n^*
	\end{align*} 
as the difference of the integrands of the last two quantities. From \eqref{ineq:1}, it follows that
	\begin{align}\label{ineq}
\int_{I_n^*} Q(\bt) d\bt \geq 0
\end{align}
	 for any  $n$ and for any boundary curve $g$. \newline
%
%
	
\noindent Define now the three sets $Q^+$, $Q^-$ and $Q^0$ as 
	\begin{align*}
	Q^+:=\{ \bt \in I_n^* :  Q(\bt)> 0\},\ Q^-:=\{ \bt \in I_n^* :  Q(\bt)< 0 \} \ \ \ \text{and} \ \ \ Q^0:=\{ \bt \in I_n^* :  Q(\bt)=0 \}.
	\end{align*}
First note that
	\[
	\lambda(Q^0)=0	 \ \ \ \text{or} \ \ \ Q\equiv 0
	\]
where $\lambda$ is the Lebesgue measure.
In the latter case, Proposition \ref{prop:adapJirak} is trivially true since then from \eqref{infDelta} we have 
	\begin{align*}
|\hat{g}(\bx)-g(\bx)| = |\hat b_{\textbf 0}-b_{\textbf 0}| = \bigtriangledown_{n}^* \leq  c_g c_1(\beta,q) h_n^\beta.
	\end{align*}
For $Q\not\equiv 0$ we have
\begin{align}\label{inf}
	\lambda(Q^+)>0	 
	\end{align}
which we will prove via a contradiction. If $\lambda(Q^+)=0$ this would imply that $\lambda(Q^-)=\lambda(I_n^*)$ since $\lambda(Q^0)=0$ and $Q^- \cup Q^+ \cup Q^0= I_n^* $. But $\frac{\lambda(Q^-)}{\lambda(I_n^*)}=1$ is a contradiction to \eqref{ineq} so \eqref{inf} must be true. This is needed for the application of Lemma \ref{inf-versus-sup}.


\noindent Note that by definition 
		\begin{align*}
	Q(\bx)=\sum_{|\bj|\leq\beta^*}(b_\bj-\hat{b}_\bj)(\textbf{0})^\bj+ \bigtriangledown_{n}^*=b_\textbf{0}-\hat{b}_\textbf{0}+ \bigtriangledown_{n}^*
	\end{align*}
	and thus
	\begin{align}\label{supQ}
	\sup_{\bt \in I_n^*}|Q(\bt)| \geq |Q(\bx)|\geq |\hat{b}_\textbf{0}-b_\textbf{0}|- \bigtriangledown_{n}^*.
	\end{align} 
	
\noindent Now define 
\beq
\tilde{Q}(\bt)=Q&\Big(&(1-t^{(1)})\max((x^{(1)}-h_n),0)+t^{(1)}\min((x^{(1)}+h_n),1),\\
&&\ldots,(1-t^{(q)})\max((x^{(q)}-h_n),0)+t^{(q)}\min((x^{(q)}+h_n),1)\Big).
\eeq
It is a polynomial on $[0,1]^q$ and it inherits from $Q$ the properties of a nonnegative integral $\int_{[0,1]^q} \tilde Q(\bt)d\bt\geq 0$ and an area with positive Lebesgue measure where $\tilde{Q}$ is positive. Thus Lemma \ref{inf-versus-sup} can be applied.
With this and \eqref{supQ} we get 
	\begin{eqnarray}\label{ineq:5}
	\nonumber c_3(\beta^*,q)(|\hat{b}_\textbf{0}-b_\textbf{0}|- \bigtriangledown_{n}^*)&\leq&c_3(\beta^*,q)\sup_{\bt \in I_n^*}|Q(\bt)|\\
	\nonumber &=&c_3(\beta^*,q)\sup_{\bt\in[0,1]^q}|\tilde{Q}(\bt)|\\
	 \nonumber &\leq&\inf_{\bt\in B_{\gamma}}\tilde{Q}(\bt)\\
	 &=&\inf_{\bt\in B_{\delta}}Q(\bt)
	\end{eqnarray}
for some constant $c_3(\beta^*,q)>0$ where $B_{\gamma}:=[\gamma_1^{(1)},\gamma_2^{(1)}]\times\cdots \times[\gamma_1^{(q)},\gamma_2^{(q)}]$ and $B_{\delta}:=[\delta_1^{(1)},\delta_2^{(1)}]\times\cdots \times[\delta_1^{(q)},\delta_2^{(q)}]$ with $\delta_1^{(r)}=(1-\gamma_1^{(r)})\max((x^{(r)}-h_n),0)+\gamma_1^{(r)}\min((x^{(r)}+h_n),1)$ and $\delta_2^{(r)}=(1-\gamma_2^{(r)})\max((x^{(r)}-h_n),0)+\gamma_2^{(r)}\min((x^{(r)}+h_n),1)$ for some ${\bf\gamma}_1,{\bf\gamma}_2\in[0,1]^q$. Thus $B_\delta$ has a volume of minimum size $c_2(\beta^*,q)h_n^q$ for some $c_2(\beta^*,q)>0$.
On the other hand, for any design point $\bX_i\in I_n^*$, by definition, 
	\begin{align}\label{ineq:6}
	Q(\bX_i) \nonumber & =\sum_{|\bj|\leq\beta^*}(b_\bj-\hat{b}_\bj)(\bX_i-\bx)^\bj+ \bigtriangledown_{n}^*\\\nonumber
	&  =g(\bX_i)-\bigtriangledown_{i}+\bigtriangledown_{n}^*-\sum_{|\bj|\leq\beta^*}\hat{b}_\bj(\bX_i-\bx)^\bj\\ \nonumber 	
	& \leq  g(\bX_i)+2 \bigtriangledown_{n}^*-Y_i\\ 
	& =|\varepsilon_i|+2 \bigtriangledown_{n}^* 
	\end{align}
	where the last inequality comes from the constraint $Y_i \leq p(\bX_i)$ and $|\bigtriangledown_{i}| \leq \bigtriangledown_{n}^*$ for all $i=1, \ldots,n$.\newline

	\noindent Combining inequalities \eqref{ineq:5} and \eqref{ineq:6}, it follows that 
	\begin{align}
\nonumber	c_3(\beta^*,q)|\hat{b}_\textbf{0}-b_\textbf{0}| & \leq c_3(\beta^*,q)\bigtriangledown_{n}^* + \inf_{\bt \in B_{\delta}}Q(\bt)\\ \nonumber
	& \leq c_3(\beta^*,q)\bigtriangledown_{n}^* + \min_{\bX_i \in B_{\delta} }Q(\bX_i)\\ 
	& \leq c_3(\beta^*,q)\bigtriangledown_{n}^* + \min_{\bX_i \in B_{\delta} }(| \varepsilon_{i}| + 2 \bigtriangledown_{n}^*). 
	\end{align}
	Then, there exists some positive constant $ C_{\beta ^*,q}$ such that 
	$$|\hat{g}(\bx)-g(\bx)|=|\hat{b}_\textbf{0}-b_\textbf{0}| \leq C_{\beta ^*,q} \bigtriangledown_{n}^* +  C_{\beta ^*,q} \min_{ \substack{ i \in \{1,\ldots,n\} \\ \bX_i \in B_{\delta} } } |\varepsilon_{i} |.$$
	
	Using the upper bound of $\bigtriangledown_{n}^*$ given in \eqref{infDelta}, there exists some positive constant $C_{\beta^*,q,c^*}$ depending only on $\beta^*$, $q$ and $c^*$ such that 
	$$|\hat{g}(\bx)-g(\bx)|=|\hat{b}_\textbf{0}-b_\textbf{0}| \leq C_{\beta^*,q,c^*}h_n^\beta + C_{\beta ^*,q} \min_{ \substack{ i \in \{1,\ldots,n\} \\ \bX_i \in B_{\delta} } } |\varepsilon_{i}| .$$
	
	Choosing a constant $J_{\beta^*, q} \in \mathbb{N}$ depending only on $q$ and $\beta^*$  large enough 
	then there exists some $\mathbf{l}=(l_1,\ldots,l_q)$ with $l_r\in\{1, \ldots, 2J_{\beta^*, q}\}$ for $r=1,\ldots,q$ such that $\bx+h_nI_{\mathbf{l}}:=\{(x^{(1)}+h_n I_{l_1})\cap[0,1]\}\times\cdots\times \{(x^{(q)}+h_nI_{l_q})\cap[0,1]\} \subseteq B_{\delta}$ where $I_k=[-1+(k-1)/J_{\beta^*, q}, -1 + k/J_{\beta^*, q}], k \geq 0$ and thus
		\begin{align*}
	|\hat{g}(\bx)-g(\bx)|\leq  C_{\beta^*,q,c^*}h_n^\beta +  C_{\beta ^*,q}  \max_{\bj\in\{1,\ldots ,2 J_{\beta^*, q}\}^q}  \left(    
	\min_{\bX_i \in \bx+h_n I_{\bj}}	|\varepsilon_{i}|   \right),
	\end{align*}
	which is the desired result and concludes the proof.
\boxi


\begin{lemma}\label{inf-versus-sup}
Let $\beta^*$ and $q$ be natural numbers. There exist positive real numbers $\delta, c$ such that for all polynomial functions $P\colon [0,1]^q\to \mathbb R$ of degree $\beta^*$ with non-negative integral over $[0,1]^q$ there exists a $\delta$-ball $B_{\delta}\subseteq [0,1]^q$ w.r.t.\ the maximum norm such that $P\geq 0$ on $B_\delta$ and
\[
    \inf_{B_\delta}P \geq c\cdot \sup_{[0,1]^q}|P|.
\]
\end{lemma}

\noindent\textbf{Proof:}
For all polynomial functions $P\colon [0,1]^q\to \mathbb R$
of degree $\beta^*$ we have
\begin{equation}\label{markov}
    ||\nabla P||_{[0,1]^q, \text{Eucl}} < 4(\beta^*)^2 \cdot \sup_{[0,1]^q}|P|,
\end{equation}
where $||\cdot||_{[0,1]^q,\text{Eucl}}$ denotes the supremum of the Euclidean norm of its argument over the set $[0,1]^q$ and $\nabla P$ stands for the gradient of $P$.
This result follows from Theorem 3.1 in Wilhelmsen (1974); see Remark \ref{rem:improv_Markov}. By the mean value theorem and the Cauchy-Schwarz inequality, this implies that for all $x,y\in[0,1]^q$ we have
\begin{equation}\label{lipschitz}
    |P(x)-P(y)| < L \cdot ||x-y|| \sup_{[0,1]^q}|P|,
\end{equation}
where $||\cdot||$ is the maximum norm and $L=4\sqrt q(\beta^*)^2$. The factor $\sqrt q$ stems from the change of norms.

Define $\varepsilon = \frac{1}{2L}, \delta = \frac{|B_\varepsilon|}{4L}$ and $c = \frac{|B_\varepsilon|}{4}$. Let $P$ be any such polynomial function. Let $P$ attain its supremum over the set $V_+ = \{x\mid P(x)>0\}$ at the point $x_0$, and denote by $x_1$ the point where $|P|$ attains its supremum on the whole cube $[0,1]^q$.
First, we show
\begin{equation}\label{sup-versus-positive-sup}
|P(x_1)| \leq \frac{2}{|B_\varepsilon|}P(x_0),
\end{equation}
where $|B_\varepsilon|$ denotes the volume of some $\varepsilon$-ball inside $[0,1]^q$. This step is necessary in case $x_1\not \in V_+$, otherwise it is not needed but the statement is still trivially true.
We note that $P$ is non-zero on the $\varepsilon$-ball $B_\varepsilon$ centered at $x_1$. Indeed, taking $y$ to be the nearest point to $x_1$ such that $P(y)=0$, Equation~\eqref{lipschitz} gives $||x_1-y|| > 2\varepsilon > \varepsilon$. Next, by the mean value theorem there exists $t_0\in B_\varepsilon$ such that
\[
    \int_{B_\varepsilon}|P(t)|\,dt = |P(t_0)||B_\varepsilon|.
\]
Applying Equation~\eqref{lipschitz} to $x_1$ and $t_0$ gives $|P(t_0)|\geq \frac{1}{2}|P(x_1)|$. Hence,
\begin{align*}
    |P(x_1)| &\leq 2|P(t_0)| \\
    &= \frac{2}{|B_\epsilon|}\int_{B_\epsilon}|P(t)|\,dt \\
    &\leq \frac{2}{|B_\epsilon|}\int_{V_+}P(t)\,dt \\
    &\leq \frac{2}{|B_\varepsilon|} P(x_0) |V_+| \\
    &\leq \frac{2}{|B_\varepsilon|} P(x_0).
\end{align*}
In the second inequality, we used that the integral of $P$ over $[0,1]^q$ is non-negative.
Next, let $B_\delta$ denote the $\delta$-ball centered at $x_0$, and let $P$ attain its infimum over $B_\delta$ at the point $x_2\in B_\delta$.
Note that $P > 0$ on $B_\delta$. Indeed, taking $y$ to be the nearest point to $x_0$ such that $P(y)=0$, Equations~\eqref{lipschitz} and  \eqref{sup-versus-positive-sup} give $||x_0-y||>2\delta > \delta$.
Finally, applying Equations~{\eqref{sup-versus-positive-sup}} and~{\eqref{lipschitz}} to $x_0$ and $x_2$ gives
\begin{align*}
    \frac{|B_\varepsilon|}{2}|P(x_1)|-P(x_2) &\leq P(x_0) - P(x_2) \\
    &\leq L |P(x_1)| \delta \\
    &= \frac{|B_\varepsilon|}{4}|P(x_1)|.
\end{align*}
Hence, $P(x_2)\geq c|P(x_1)|$ as required. This concludes the proof. \boxi

\begin{remark}\label{rem:improv_Markov}
    The bound given in the Markov inequality in Equation \eqref{markov} in Lemma \ref{inf-versus-sup} is far from being optimal. For the sake of completeness, we recall the full statement of Theorem 3.1 in Wilhelmsen (1974). With the above notation, it writes 
    \begin{align}\label{Wilhelsem}
            ||\nabla P||_{T, Eucl} \leq \frac{4 \beta^{*2}}{\omega(T)} \sup_{T}|P| 
    \end{align}
    where  $T$ is a compact and convex set with non-empty interior and $\omega(T)$ stands for the thickness of $T$, that is the minimum distance between two parallel supporting hyperplanes for $T$. Note that in the context of Lemma \ref{inf-versus-sup}, $\omega (T)=\omega ([0,1]^q)=1$. The bound of Equation \eqref{Wilhelsem} has been improved later on  by Kro\'{o}  and  R\'{e}v\'{e}sz (1999). They showed that for any convex body (that is a compact convex set with nonempty interior) the constant $\frac{4 \beta^{*2}}{\omega(T)}$ may be replaced by $\frac{4 \beta^{*2} -2\beta^{*}}{\omega(T)}$. One may go further observing that the space $[0,1]^q$, under some adequate normalisation, may be seen as a central symmetric convex body (a convex body is central symmetric if and only if with proper shift it is the unit
ball of some norm on $\mathbb{R}^q$) so one may use the result of Sarantopoulos  (1991) to get an even smaller bound. We refer to the book of Rassias and Tóth (2014) for an exhaustive treatment of Markov-type inequalities for multivariate polynomials. For the sake of simplicity, we have presented the result with the bound given in Wilhelmsen (1974).

\end{remark}
\begin{proposition}\label{min-eps-rand}
Assume that model \eqref{model-random} holds under (K2)-(K4).
Then 
\[\sup_{\bx\in[0,1]^q}\min_{\substack{i\in\{1,\ldots,n\}\\ \bX_{i}\in\bx+h_nI}}|\e_{i}|=O_{\P}\left(\left(\frac{\left|\log\left(h_n\right)\right|}{nh_n^q}\right)^{\frac 1{\alpha}}\right)\]
for every non-degenerate subinterval $I \subseteq[-1,1]^q$ where we set $\min_\emptyset|\e_{i}|:=0$.\end{proposition}

\noindent Before we prove this proposition let us first note that under assumption (K3)

$$  \frac{\left|\log\left(h_n\right)\right|}{nh_n^q} = O\left( \frac{\log(n)}{nh_n^q}\right)$$ 
so that Proposition \ref{min-eps-rand} leads to the appropriate rates in the main result of Theorem \ref{th:rateunif}. Similarly Proposition \ref{min-eps} and Proposition \ref{min-eps-nonequi} respectively gives way to the expected rates in the main result of Theorem \ref{th2:rateunif}.
\\ \\
\noindent\textbf{Proof:}
The proof is similar to the proof of Proposition \ref{min-eps}. There are just some preliminaries to consider to deal with the random case.

\noindent Observe that it is obvious that we only have to consider those $I$ where there exists a $d>0$ such that $|(x^{(r)}+h_nI_r)\cap[0,1]|\geq dh_n$ for all $r=1,\ldots,q$ with $I=I_1\times \cdots\times I_q$ and by $|\cdot|$ we mean the length of the one-dimensional intervals. 

\noindent Next we will show that with probability converging to one there are at least $d_n=O(nh_n^q)$ random design points in every hypercube $I_n\subseteq[0,1]^q$ with edge length $dh_n$. This implies that at least $d_n$ random design points lie in $\bx+h_nI$.
The number of points in $I_n$ can be written as $\sum_{i=1}^nI\{\bX_i\in I_n\}$ and for all $I_n$
\beq
\frac 1n\sum_{i=1}^nI\{\bX_i\in I_n\}&\geq&\P(\bX_1\in I_n) \\
&&-\left|\frac 1n\sum_{i=1}^nI\{\bX_i\in I_n\}-\E\left[\frac 1n\sum_{i=1}^nI\{\bX_i\in I_n\}\right]\right|.
\eeq
For the term $\P(\bX_1\in I_n)$  the lower bound $(dh_n)^q\inf_{\bt\in[0,1]^q}f_\bX(\bt)$ by 
\[\P(\bX_1\in I_n)=\int_{I_n}f_\bX(\bt)d\bt\]
and the consideration about the length of $I_n$. The bound on $\P(\bX_1\in I_n)$ is of order $O(h_n^q)$ by assumption (K4). Thus it remains to prove that
\begin{equation}\label{min-eps-rand-pollard}
\sup_{I_n}\left|\frac 1n\sum_{i=1}^nI\{\bX_i\in I_n\}-\E\left[\frac 1n\sum_{i=1}^nI\{\bX_i\in I_n\}\right]\right|=O(h_n^q)
\end{equation}
with probability converging to one. Set therefore $P_nf_{I_n}:=\frac 1n\sum_{i=1}^nf_{I_n}(\bX_i)$ and $Pf_{I_n}:=\E[f_{I_n}(\bX_1)]$ with $f_{I_n}(\bt):=I\{\bt\in I_n\}$. Note that $Pf_{I_n}^2=\P(\bX_1\in I_n)\leq (dh_n)^q\sup_{\bt\in[0,1]^q}f_\bX(\bt)$. Note further that $f_{I_n}(\bt)=I\{\bt\in I_n\}=I\{\|\bt-\frac {\textbf{a}_n+\textbf{b}_n}2\|\leq \frac{d}2h_n\}$  - where $\textbf{a}_n=(a_{n,1},\ldots,a_{n,q})$, $\textbf{b}_n=(b_{n,1},\ldots,b_{n,q})$ and $I_n=[a_{n,1},b_{n,1}]\times\cdots\times[a_{n,q},b_{n,q}]$ - and thus the conditions of Example 38 and Problem 28 in Pollard (1984) are fullfilled. Then since $|f_{I_n}|\leq 1$ Theorem 37 in Pollard (1984) can be applied and thus
\[\sup_{I_n}|P_nf_{I_n}-Pf_{I_n}|=o(h_n^q)\qquad\text{a.\,s.}\]
which proves \eqref{min-eps-rand-pollard}.
Now the assertion of the Proposition follows with the same arguments as in the proof of Proposition \ref{min-eps} with the modification from the proof of Proposition \ref{min-eps-nonequi}. \boxi

\subsection{Proof of Theorem \ref{th2:rateunif}}
The proof of Theorem \ref{th2:rateunif} is similar to the proof of Theorem \ref{th:rateunif} and is based on the following Propositions \ref{prop2:adapJirak}, \ref{min-eps} and \ref{min-eps-nonequi} respectively.
\begin{proposition}\label{prop2:adapJirak}
	Assume that model \eqref{model0} holds under (K2) and (K4') and let $\hat{g}$ be defined in \eqref{estimatorIntegralFix}-\eqref{eq:integFix} satisfying (K1) for some $\beta \in (0, \beta^* +1]$ and some $c_g \in [0,c^*]$. Then, there exists constants $C_{\beta^*,q,c^*},C_{\beta^*,q}$ and a natural number $J_{\beta^*, q}$ such that  for all $\bx \in [0,1]^q $
	\begin{align*}
|\hat{g}(\bx)-g(\bx)|\leq C_{\beta^*,q,c^*}h_n^\beta +  C_{\beta ^*,q}  \max_{\bj\in\{1,\ldots ,2 J_{\beta^*, q}\}^q}  \left(    
	\min_{\substack{i\in\{1,\ldots,n\}\\ \bx_{i,n} \in \bx+h_n I_{\bj}}}	|\varepsilon_{i,n}|   \right),
	\end{align*}
	where $\bx+h_nI_\bj=\{(x^{(1)}+h_n I_{j_1})\cap[0,1]\}\times\cdots\times \{(x^{(q)}+h_nI_{j_q})\cap[0,1]\}$ with $I_k=[-1+(k-1)/J_{\beta^*, q}, -1 + k/J_{\beta^*, q}]$.
\end{proposition}

\noindent\textbf{Proof:}
	The proof is almost identical to the proof of Proposition \ref{prop:adapJirak} and is skipped here for the sake of conciseness.

\begin{proposition}\label{min-eps}
Assume that model \eqref{model0} holds with equidistant design points under (K2) and (K3') (which is equivalent to (K3) in this case).
Then 
\[\sup_{\bx\in[0,1]^q}\min_{\substack{i\in\{1,\ldots,n\}\\ \bx_{i,n}\in\bx+h_nI}}|\e_{i,n}|=O_{\P}\left(\left(\frac{\left|\log\left(h_n\right)\right|}{nh_n^q}\right)^{\frac 1{\alpha}}\right)\]
for every non-degenerate subinterval $I \subseteq[-1,1]^q$ where we set $\min_\emptyset|\e_{i,n}|:=0$.
\end{proposition}

\noindent\textbf{Proof:} The proof is similar to the proof of Lemma A.2 in Drees et al.\ (2019) and the proof of Lemma A.1 in Neumeyer et al.\ (2019) but comprehensive adaptions are needed to deal with the multivariate case.

\noindent Let $Z_1,Z_2,\ldots$ be iid with the same distribution as $-\e_{i,n}$ with cumulative distribution function $U$. Recall that the distribution of the errors $\e_{i,n}$ do not depend of $n$, that is why, for the sake of clarity, the second index $n$ is omitted in the $Z_i$. To prove the result we shall show that 
\[\exists L<\infty:\quad \P\left(\sup_{\bx\in[0,1]^q}\min_{\substack{i\in\{1,\ldots,n\}\\ \bx_{i,n}\in\bx+h_nI}}Z_i>Lr_n\right)\nto 0\]
with $r_n=\left(\frac{\left|\log\left(h_n\right)\right|}{nh_n^q}\right)^{\frac 1{\alpha}}$.

\noindent It is obvious that we only have to consider those $I$ where there exists a $d>0$ such that $|(x^{(r)}+h_nI_r)\cap[0,1]|\geq dh_n$ for all $r=1,\ldots,q$ with $I=I_1\times \cdots\times I_q$ and by $|\cdot|$ we mean the length of the one-dimensional intervals. Note that at least $d_n=\lfloor dh_nn^{\frac 1q}\rfloor^q$ design points lie in such an $\bx+h_nI$.

\noindent Now we want to arrange the design points $\bx_{i,n}$ in sets of size smaller than $d_n$ such that we can replace $\sup_{\bx\in[0,1]^q}$ by a maximum over these sets. To this end
for each $n$ we choose $H_n\in(\frac {dh_n}4,\frac{dh_n}2]$ such that $\frac 1{H_n}$ is an integer. With this we build
$$D_{\bl,n}:=\{\bx_{i,n}:x_{i,n}^{(r)}\in((l^{(r)}-1)H_n,l^{(r)}H_n], r=1,\ldots,q\}$$
for all $\bl\in\{1,\ldots,\frac 1{H_n}\}^q$.
Note that these sets contain $d_{H_n}:=\lfloor H_n n^{\frac 1q}\rfloor ^q< d_n$ design points and that for every $\bx \in [0,1]^q$ there exists $D_{\bl,n}$ with $D_{\bl,n}\subseteq\{\bx_{i,n}: \bx_{i,n}\in\bx+ h_nI\}$.

\noindent Now the supremum can be replaced by maximum, i.\ e.\ for all $y>0$
\beq
\P\left(\sup_{\bx\in[0,1]^q}\min_{\substack{i\in\{1,\ldots,n\}\\ \bx_{i,n}\in\bx+ h_nI}}Z_i>y\right)&\leq& \P\left(\max_{\bl\in\{1,\ldots,\frac 1{H_n}\}^q}\min_{ \substack{i\in\{1,\ldots,n\}\\ \bx_{i,n}\in D_{\bl,n}}}Z_i>y\right)\\
&=& 1-\P\left(\min_{ \substack{i\in\{1,\ldots,n\}\\ \bx_{i,n}\in D_{\bl,n}}}Z_i\leq y\ \forall \bl\in\left\{1,\ldots,\frac 1{H_n}\right\}^q\right)\\
&=&1-\prod_{\bl\in\left\{1,\ldots,\frac 1{H_n}\right\}^q}\P\left(\min_{ \substack{i\in\{1,\ldots,n\}\\ \bx_{i,n}\in D_{\bl,n}}}Z_i\leq y\right)
\eeq
since the sets $D_{\bl,n}$ are disjoint and thus the corresponding $Z_i$ are independent. Further $\P\left(\min_{ \substack{i\in\{1,\ldots,n\}\\ \bx_{i,n}\in D_{\bl,n}}}Z_i> y\right)=\overline{U}(y)^{d_{H_n}}$ and thus
\[\P\left(\sup_{\bx\in[0,1]^q}\min_{\substack{i\in\{1,\ldots,n\}\\ \bx_{i,n}\in\bx+ h_nI}}Z_i>y\right)\leq 1-(1-\overline{U}(y)^{d_{H_n}})^{\frac 1{H_n^q}}.\]
It remains to show that for sufficiently large $L$
\begin{equation}\label{min-eps2}
1-\left(1-\overline{U}(Lr_n)^{d_{H_n}}\right)^{\frac 1{H_n^q}}\nto 0
\end{equation}
which is true if 
\begin{equation}\label{min-eps1}
\overline{U}(Lr_n)^{d_{H_n}}=o\left(H_n^q\right).
\end{equation}
This implication follows from several Taylor expansions. Let $z_n\to 0$ and $y_n\to\infty$ denote some sequences with $y_nz_n\to0$. Now by a second order Taylor expansion of $\log(1-z_n)$ we get for $n$ sufficiently large
\beq
1-(1-z_n)^{y_n}&=&1-\exp\left(y_n\log\left(1-z_n\right)\right)\\
&=&1-\exp(-y_nz_n)+O(y_nz_n^2)\\
&=&1-\exp(0)+O(y_nz_n)+O(y_nz_n^2)\\
&=&O(y_nz_n)
\eeq
where the third equality follows by a first order Taylor expansion of $\exp(-y_nz_n)$. This means that
$$1-\left(1-\overline{U}(Lr_n)^{d_{H_n}}\right)^{\frac 1{H_n^q}}=O\left(\frac 1{H_n^q}\overline{U}(Lr_n)^{d_{H_n}}\right)$$
which proves our claim $(\ref{min-eps1})\Rightarrow (\ref{min-eps2})$.

\noindent Note that $U(Lr_n)=cL^{\alpha}r_n^{\alpha}(1+o(1))$ for some positive constant $c$ and thus we have by a second order Taylor expansion of $\log(1-U(Lr_n))$
\beq
\overline{U}(Lr_n)^{d_{H_n}}&=&\exp(d_{H_n}\log(1-U(Lr_n)))\\
&=&\exp\left(-cL^{\alpha}d_{H_n}r_n^{\alpha}(1+o(1))\right)\\
&=&\exp\left(-cL^{\alpha}\frac{d_{H_n}}{nh_n^q}|\log(h_n)|(1+o(1))\right)\\
&\leq&\exp\left(-(q+\delta)|\log(h_n)|\right)\\
&=&o(h_n^q)
\eeq
for some $\delta>0$ and sufficiently large $L$ and $n$, where the second to last line follows from the fact that one can find a $\delta >0$ such that $cL^{\alpha}\frac{d_{H_n}}{nh_n^q}(1+o(1))\geq q+\delta$ provided that both $L$ and $n$ are sufficiently large. This concludes the proof since $h_n=O(H_n)$.

\boxi

\begin{proposition}\label{min-eps-nonequi}
Assume that model \eqref{model0} holds under (K2), (K3') and (K4').
Then 
\[\sup_{\bx\in[0,1]^q}\min_{\substack{i\in\{1,\ldots,n\}\\ \bx_{i,n}\in\bx+h_nI}}|\e_{i,n}|=O_{\P}\left(\left(\frac{\left|\log\left(h_n\right)\right|}{d_n}\right)^{\frac 1{\alpha}}\right)\]
with $d_n=d_n(1)$, for every non-degenerate subinterval $I \subseteq[-1,1]^q$ where we set $\min_\emptyset|\e_{i,n}|:=0$.
\end{proposition}

\noindent\textbf{Proof:}
As in the proof of Proposition \ref{min-eps} it is obvious that we only have to consider those $I$ where there exists a $d>0$ such that $|(x^{(r)}+h_nI_r)\cap[0,1]|\geq dh_n$ for all $r=1,\ldots,q$. By assumption (K4') there lie at least $d_n(d)$ design points in such an $\bx+h_nI$. 

\noindent The proof is similar to the proof of Proposition \ref{min-eps}. We only have to adjust some steps of the proof to the non-equidistant case. Note to this end that the number of design points in the sets $D_{\bl,n}$ from the proof of Proposition \ref{min-eps}  may differ from set to set in the non-equidistant case but that there are at least $d_{H_n}:=d_n(H_n/h_n)$  points in each set by assumption (K4').

\noindent For all $\bl\in\{1,\ldots,\frac 1{H_n}\}^q$ it holds $\P\left(\min_{ \substack{i\in\{1,\ldots,n\}\\ \bx_{i,n}\in D_{\bl,n}}}Z_i> y\right)=\overline{U}(y)^{\#D_{\bl,n}}\leq\overline{U}(y)^{d_{H_n}}$ where $\#D_{\bl,n}$ is the number of design points in $D_{\bl,n}$.
Thus like in the proof of Proposition \ref{min-eps} we have
\[\P\left(\sup_{\bx\in[0,1]^q}\min_{\substack{i\in\{1,\ldots,n\}\\ \bx_{i,n}\in\bx+ h_nI}}Z_i>y\right)\leq 1-(1-\overline{U}(y)^{d_{H_n}})^{\frac 1{H_n^q}}.\]
Now the assertion follows with the same arguments as in the proof of Proposition \ref{min-eps}.
\boxi

\section*{Acknowledgements}
We would like to thank Holger Drees, Natalie Neumeyer, Bernd Sturmfels, Lorenzo Venturello, Avinash Kulkarni and János Kollár for helpful discussions.
Further we would like to thank the anonymous referees for many constructive remarks that led to a significant improvement of the paper. Financial support by the DFG (Research Unit FOR 1735 {\it Structural Inference in Statistics: Adaptation and Efficiency}) is gratefully acknowledged.

\section*{References}

\begin{description}

\item Aigner, D., Lovell, C. K. and Schmidt, P. (1977). Formulation and estimation of stochastic frontier production function models. \textit{Journal of Econometrics}, \textbf{6}(1), 21-37.

\item Andersen, E. D. and Andersen, K. D. (2000). The MOSEK interior point optimizer for linear programming: an implementation of the homogeneous algorithm. \textit{High Performance Optimization}, 197--232

\item Box, G. E. P. and Cox, D. R. (1964). An analysis of transformations. \textit{Journal of the Royal  Statistical Society Series B},  \textbf{26}, 211--252.

\item Brown, L.D. and Low M.G. (1996). Asymptotic equivalence of nonparametric regression and white noise. \textit{The Annals of Statistics}, \textbf{24}, 2384--2398.

 \item Cooper, W., Seiford, L. M. and Zhu, J. (Eds.) (2011). \textit{Handbook on Data Envelopment Analysis} (Vol. 164). Springer Science and Business Media.
 
 \item Cornwell, C. and Schmidt, P. (2008). Stochastic frontier analysis and efficiency estimation. \textit{In The Econometrics of Panel Data} 697-726. Springer, Berlin, Heidelberg.
 
 \item Daouia, A. and Simar, L. (2005). Robust nonparametric estimators of monotone boundaries. \textit{Journal of Multivariate Analysis}, \textbf{96}(2), 311-331.

\item Daouia, A., Noh, H. and Park, B. U. (2016).  Data envelope fitting with constrained polynomial splines.
\textit{Journal of the Royal Statistical Society, Series B} \textbf{78}, 3--30.

\item Deprins, D., Simar, L. and Tulkens, H. (1984). \textit{Measuring labor-efficiency in post offices}. No 571, CORE Discussion Papers RP, Université catholique de Louvain, Center for Operations Research and Econometrics (CORE).

 \item de Haan, L. and Resnick, S. (1994). Estimating the home range. \textit{Journal of Applied
Probability}, \textbf{31}, 700--720.

\item de Haan, L. and Ferreira, A. (2006). \textit{Extreme Value Theory: an Introduction}. Springer
Series in Operations Research and Financial Engineering. Springer, New York.

 \item Drees, H., Neumeyer, N. and Selk, L. (2019). Estimation and hypotheses tests in boundary regression models. \textit{Bernoulli}, {\bf 25}, 424--463.
 
 \item Farrell, M. J. (1957). The measurement of productive efficiency. \textit{Journal of the Royal Statistical Society: Series A (General)}, \textbf{120}(3), 253-281.
 
 \item Francis, J. C. and Kim, D. (2013). \textit{Modern Portfolio Theory: Foundations, Analysis, and New Developments} (Vol. 795). John Wiley \& Sons.

 \item Gardes, L. (2002). Estimating the support of a poisson process via the faber-schauder basis and extreme values. In \textit{Annales de l'ISUP} \textbf{46}(1-2), 43--72.
 
  \item Gijbels, I., Mammen, E., Park, B.U. and Simar, L. (1999). On estimation of monotone
and concave frontier functions. \textit{Journal of American Statistical Association}, \textbf{94}, 220--228.

\item Gijbels, I. and Peng, L. (2000). Estimation of a support curve via order statistics. \textit{Extremes}, \textbf{3}, 251--277.
 
 \item Girard, S. and Jacob, P. (2003). Projection estimates of point processes boundaries. \textit{Journal of Statistical Planning and Inference}, \textbf{116}(1), 1--15.
 
 \item  Girard, S. and Jacob, P. (2004). Extreme values and kernel estimates of point processes boundaries. \textit{ESAIM: Probability and Statistics}, \textbf{8}, 150--168.
 
 \item Girard, S., Iouditski, A. and Nazin, A. V. (2005). L1-optimal nonparametric frontier estimation via linear programming. \textit{Automation and Remote Control},  \textbf{66}(12), 2000--2018.

  \item Girard, S. and Jacob, P. (2008). Frontier estimation via kernel regression on high powertransformed data. \textit{Journal of Multivariate Analysis} \textbf{99}, 403–420.
 
  \item Girard, S., Guillou, A. and Stupfler, G. (2013). Frontier estimation with kernel regression
on high order moments. \textit{Journal of Multivariate Analysis}, \textbf{116}, 172–189.

\item Goetzmann, W. N., Brown, S. J., Gruber, M. J. and Elton, E. J. (2014). \textit{Modern Portfolio Theory and Investment Analysis}. John Wiley \& Sons, 237.
 
 \item Hall, P., Nussbaum, M. and Stern, S.E. (1997). On the estimation of a support curve
of indeterminate sharpness. \textit{Journal of Multivariate Analysis}, \textbf{62}, 204–232.

 \item Hall, P., Park, B.U. and  Stern, S.E. (1998). On polynomial estimators of frontiers and boundaries. \textit{Journal of Multivariate Analysis},  \textbf{66}, 71--98.

\item Hall, P. and Park, B.U. (2004). Bandwidth choice for local polynomial estimation of
smooth boundaries. \textit{Journal of Multivariate Analysis}, \textbf{91} (2), 240–-261.
 
 \item Hall, P. and Van Keilegom, I. (2009). Nonparametric ``regression'' when errors are positioned at end-points. \textit{Bernoulli}, \textbf{15}, 614--633.
  
 \item H\"ardle, W., Park, B.U. and Tsybakov, A.B. (1995). Estimation of non-sharp support
boundaries. \textit{Journal of Multivariate Analysis}, \textbf{55}, 205–-218.

\item Jacob, P. and Suquet, C. (1995). Estimating the edge of a Poisson process by orthogonal series. \textit{Journal of Statistical Planning and Inference}, \textbf{46}(2), 215--234.
 
 \item Jirak, M., Meister, A.\ and Rei{\ss}, M. (2014). Adaptive estimation in nonparametric regression with one-sided errors. \textit{The Annals of Statistics}, {\bf 42}(5), 1970--2002.
 
 \item Jones, M.C. and Pewsey, A. (2009). Sinh-arcsinh distributions. \textit{Biometrika},  \textbf{96}, 761--780.
 
 \item Kelly, E., Shalloo, L., Geary, U., Kinsella, A. and  Wallace, M. (2012). Application of data envelopment analysis to measure technical efficiency on a sample of Irish dairy farms. \textit{Irish Journal of Agricultural and Food Research}, 63-77.
 
 \item Knight, K. (2001). Limiting distributions of linear programming estimators. \textit{Extremes}, \textbf{4} (2), 87--103.

\item Korostelëv A.P. and Tsybakov A.B. (1993). \textit{Minimax Theory of Image Reconstruction}. Lecture Notes in Statistics,  \textbf{82}. Springer, New York, NY

\item Korostelëv, A. P., Simar, L. and  Tsybakov, A. B. (1995). Efficient estimation of monotone boundaries. \textit{The Annals of Statistics}, \textbf{23}(2), 476-489.
 
 \item Kro\'{o}, A. and  R\'{e}v\'{e}sz, S. (1999). On Bernstein and Markov-type inequalities for multivariate polynomials on convex bodies. \textit{Journal of approximation theory}, \textbf{99}(1), 134-152.
 
  \item Kumbhakar, S. C. and Lovell, C. K. (2003). \textit{Stochastic Frontier Analysis}. Cambridge University Press.
  
  \item Kumbhakar, S. C., Park, B. U., Simar, L. and Tsionas, E. G. (2007). Nonparametric stochastic frontiers: a local maximum likelihood approach. \textit{Journal of Econometrics}, \textbf{137}(1), 1--27.
  
  \item Linton, O., Sperlich, S. and Van Keilegom, I. (2008). Estimation on a semiparametric transformation model. \textit{Annals of  Statisics}, \textbf{36}, 686–718.
 
 \item Lovell, C. K., Grosskopf, S., Ley, E., Pastor, J. T., Prior, D. and Eeckaut, P. V. (1994). Linear programming approaches to the measurement and analysis of productive efficiency. \textit{Top}, \textbf{2}(2), 175--248.
 
 \item Mardani, A., Streimikiene, D., Balezentis, T., Saman, M. Z. M., Nor, K. M. and Khoshnava, S. M. (2018). Data envelopment analysis in energy and environmental economics: An overview of the state-of-the-art and recent development trends. \textit{Energies} 2018, \textbf{11}(8).
 
  \item Meister, A.\ and Rei{\ss}, M. (2013). Asymptotic equivalence for nonparametric regression with non-regular errors. \textit{Probability Theory and Related Fields}, {\bf 155}, 201--229.
 
 \item Menneteau, L. (2008). Multidimensional limit theorems for smoothed extreme value estimates of point processes boundaries. \textit{ESAIM: Probability and Statistics}, \textbf{12}, 273-307.

\item Meeusen, W. and van Den Broeck, J. (1977). Efficiency estimation from Cobb-Douglas production functions with composed error. \textit{International Economic Review}, 435-444.
 
 \item M\"{u}ller, U.U.\ and Wefelmeyer W. (2010). Estimation in nonparametric regression with non-regular errors. \textit{Communications in Statistics - Theory and  Methods}, {\bf 39}, 1619--1629.
 
 \item Narimatsu, H., Nakata, Y., Nakamura, S., Sato, H., Sho, R. et al. (2015). Applying data envelopment analysis to preventive medicine: a novel method for constructing a personalized risk model of obesity. \textit{PLoS One}, \textbf{10}(5), e0126443.
 
 \item Neumeyer, N., Selk, L. and Tillier, C. (2019). Semi-parametric transformation boundary regression models. \textit{Annals of the Institute of Statistical Mathematics}, \textbf{27}, 1--29.
 
 \item Pollard, D. (1984). \textit{Convergence of Stochastic Processes}. Springer, New York.
 
 \item Ramanathan, R. (2003). \textit{An Introduction to Data Envelopment Analysis: a Tool for Performance Measurement}. Sage.
 
\item  Rassias, T. M. and T\'{o}th, L. (Eds.). (2014). \textit{Topics in Mathematical Analysis and Applications}. Springer International Publishing. 

\item Reiß, M. and  Selk, L. (2017). Efficient estimation of functionals in nonparametric boundary models. \textit{Bernoulli}, \textbf{23}(2), 1022--1055.
 
 \item Sarantopoulos, Y. (1991). Bounds on the derivatives of polynomials on Banach spaces. \textit{Mathematical Proceedings of the Cambridge Philosophical Society}, \textbf{110}(2), 307--312. 
 
 \item Simar, L. and Wilson, P.W. (1998). Sensitivity analysis of efficiency scores: how to bootstrap in nonparametric frontier models. \textit{Management Science}, \textbf{44}, 49--61.
 
 \item Simar, L. and Zelenyuk, V. (2011). Stochastic FDH/DEA estimators for frontier analysis. \textit{Journal of Productivity Analysis}, \textbf{36}(1), 1--20.
 
 \item Tsybakov, A. B. (1994). Multidimensional change-point problems and boundary estimation. \textit{Change-Point Problems, IMS Lecture Notes} \textbf{23},
 
\item Vaidya, P. M. (1989). Speeding-up linear programming using fast matrix multiplication. \textit{30th Annual Symposium on Foundations of Computer Science}, 332--337.

\item Wilhelmsen, D. R. (1974). A Markov inequality in several dimensions. \textit{Journal of Approximation Theory}, \textbf{11}(3), 216--220.

\item Wilson, P.W. (2003). Testing independence in models of productive efficiency. \textit{Journal of Productivity Analysis},  \textbf{20}, 361--390.

\end{description}

\noindent Corresponding author:\\ Leonie  Selk, University of Hamburg, Department of Mathematics, Bundesstrasse 55, 20146 Hamburg, Germany.  \texttt{Email: leonie.selk@math.uni-hamburg.de }

\end{document}